\documentclass[12pt,leqno]{article}

\usepackage{amssymb}
\usepackage{amsmath,amsthm}
\usepackage{latexsym}

\renewcommand{\span}{\mbox{\rm span}}
\newcommand{\vv}{\vec{v}}
\newcommand{\strip}{{\mbox{\rm strip}}}

\newcommand{\nn}{\nonumber}
\newcommand{\FT}{{\mathcal F}}
\newcommand{\la}{\langle}
\newcommand{\ra}{\rangle}
\newcommand{\diam}{\mbox{\rm diam}}
\newcommand{\sgn}{\mbox{\rm sgn}}
\newcommand{\card}{\mbox{\rm card}}
\newcommand{\bea}{\begin{eqnarray}}
\newcommand{\enea}{\end{eqnarray}}
\newcommand{\be}{\begin{equation}}
\newcommand{\ee}{\end{equation}}
\newcommand{\supp}{\mbox{\rm supp}}
\newcommand{\eproof}{\hspace{.1in}\rule{.1in}{.1in}\\ \smallskip}
\renewcommand{\to}{\rightarrow}
\newcommand{\half}{\frac{1}{2}}

\newcommand{\Cb}{\overline{C}}

\newcommand{\plates}{{\mathcal P}}
\newcommand{\lines}{{\mathcal L}}

\newcommand{\eps}{\varepsilon}

\newcommand{\R}{{\mathbb R}}
\newcommand{\Z}{{\mathbb Z}}
\newcommand{\Prob}{{\mathbb P}}
\newcommand{\A}{{\mathcal A}}
\newcommand{\B}{{\mathcal B}}
\newcommand{\C}{{\mathcal C}}
\newcommand{\W}{{\mathcal W}}
\newcommand{\El}{{\mathcal E}}

\newcommand{\CCet}{\C^C_{\epsilon t}}

\newcommand{\Acet}{\A^C_{\epsilon t}}

\newcommand{\durch}{\C^{C_1}_{\epsilon t}\cap\C^{C_2}_{\epsilon t}}

\newcommand{\proj}{\mbox{\rm proj}}
\newcommand{\dist}{\mbox{\rm dist}}
\newcommand{\good}{{\mathcal G}}
\newcommand{\Erw}{{\mathbb E}}
\newcommand{\les}{\lesssim}

\newcommand{\wt}{\widetilde}

\newtheorem{theorem}{Theorem}
\newtheorem{lemma}{Lemma}

\newtheorem{cor}{Corollary}
\newtheorem{defi}{Definition}

\begin{document}

\begin{center}
{\LARGE{On continuum incidence problems related to\\  harmonic analysis\\}}
{\mbox{      }\\}
{\large{W. Schlag}}
\end{center}

\section{Introduction}

The purpose of this note is to put some recent work in
harmonic analysis involving combinatorics of circles and spheres in perspective.
A review of some of this work was given by Wolff in~\cite{tomrev}, 
and we will not duplicate what is said there. On the other hand,
Wolff's recent paper~\cite{W2} allows for some simplifications beyond~\cite{tomrev}.
More precisely, the difficult argument in~\cite{W} can be considerably
simplified using an inequality from~\cite{W2}. We present a detailed argument in Section~\ref{sec:L3} below.
Furthermore, in order to illustrate the use of combinatorial arguments 
dealing with circles (or spheres), we present some simple cases where these
arguments are very transparent, but do not appear to be known. 
One such case is the standard estimate
\begin{equation}
\label{eq:3/2} 
\|\sigma \ast f\|_{L^3(\R^2)} \lesssim \|f\|_{L^{\frac32}(\R^2)}, 
\end{equation}
where $\sigma$ is the surface measure on the circle, and similarly in higher
dimensions. We show how to obtain restricted weak-type versions of these
bounds in a very direct way without the Fourier transform. Another proof of this fact
that does not rely on the Fourier transform was found by Oberlin~\cite{O1}, 
using a multilinear interpolation scheme that originates in Christ~\cite{christ}.
The argument that is used in this paper is rather different, and it also 
extends easily to families of surfaces satisfying the rotational curvature condition of Phong and Stein, see~\cite{PhSt}.
The new feature here is that only two derivatives are needed on the defining function of the family
of surfaces (in all dimensions), whereas all previous methods require 
a large number of derivatives that goes to infinity with the dimension of the ambient space.
Another case we discuss is the classical Strichartz estimate for the wave-equation, namely
\begin{equation}
\label{eq:26} 
\|u\|_{L^6(\R_x^2\times\R_t)} \les \|f\|_{\dot{H}^{\half}(\R^2)} 
\text{\ \ provided\ \ }\Box u=0,\;u(0,\cdot)=f,\;\partial_t u(0,\cdot)=0. 
\end{equation}
It turns out that a weaker version of the estimate can be derived without using
the Fourier transform at all, see Section~\ref{sec:strich}. There has been 
some recent interest in proofs of the Strichartz estimates that do not rely
on the Fourier transform, see Klainerman~\cite{Kl} and Smith, Tataru~\cite{tatsmith}. 
Hence it might be of some value
to give further instances where this can be done. The point of our method, based on
Marstrand's work~\cite{Mar} and Kolasa, Wolff~\cite{KW}, is to recast~\eqref{eq:26} 
as a statement about a certain multiplicity function for a large collection of
annuli in the plane and then to bound this multiplicity by means of geometric-combinatorial
arguments. One possible advantage of this approach is that it carries over to a
variable coefficient setting. More precisely, one can replace circles by curves 
satisfying Sogge's cinematic curvature condition from~\cite{sogge}.
As pointed out in~\cite{tomrev} and~\cite{W2}, some questions about tangent
circles remain open that appear to be quite difficult. In the final section of
this paper we present  examples that show why a conjecture of Wolff (and possibly others) 
would be best possible. These examples are quite standard, but it is not clear to the author
if they have appeared in print before.

\section{Averages over hypersurfaces with nonvanishing rotational curvature}

We now turn to~\eqref{eq:3/2}. 
Note that no such estimate can hold for the boundary measure of a square.
The estimate~\eqref{eq:3/2} is usually proved by means of complex interpolation, see Stein~\cite{Stein1}. 
In \cite{O1} Oberlin, however, found a proof of the restricted weak-type bound that
does not rely on the Fourier transform using an idea of multilinear interpolation of Christ~\cite{christ}.
Recently there has been a lot of activity around smoothing properties of averages
with respect to curves and surfaces. In particular, we would like to point out the
work by Christ~\cite{christ2}, Oberlin~\cite{O2}, and Bak, Oberlin, Seeger~\cite{BOS}.
Here we develop another argument that is quite different from the approach of these works. 
We first present the argument  for circles in two dimensions, and then generalize it to higher dimensions
and surfaces obeying the rotational curvature condition of Phong and Stein.

\noindent Let $C(x,r)$ be a circle centered at $x\in \R^2$ with radius $r\in[1,2]$. Denote
the $\delta$-neighborhood of $C(x,r)$ by $C^\delta(x,r)$. 
Fix any small $\delta>0$ and let $E\subset[0,1]^2$ (it suffices to consider that case). 
Fix $\lambda>0$  and define
\begin{equation}
\label{eq:Fdef} 
F=\{x\in\R^2\:|\: (\chi_E \ast \sigma_\delta)(x) > \lambda\}. 
\end{equation}
Here $\sigma_\delta$ is the normalized measure on $C^\delta(0,1)$. We need to show that
\begin{equation}
\label{eq:restwt}
|F| \le C\,\lambda^{-3}|E|^2
\end{equation}
with some absolute constant $C$. This is precisely the restricted weak-type form of~\eqref{eq:3/2}.\\ 
We first discretize $E$ on scale $\delta$. Partition $[0,1]^2$ into squares~$Q_j$ 
of side-length~$\delta$ and let 
\[ E_\ell = \bigcup_{j\::\:2^{-\ell}\delta^2 < |Q_j\cap E| \le 2^{-\ell+1}\delta^2} Q_j\cap E\]
for $\ell\ge1$. Clearly, $E=\bigcup_\ell E_\ell$ and we define 
\[ \wt{E_\ell} = \bigcup_{j\::\:2^{-\ell}\delta^2<|Q_j\cap E|\le 2^{-\ell+1}\delta^2} Q_j.\]
Note that 
\begin{equation}
\label{eq:Eellmes} 
2^{\ell}|E_\ell| < |\wt{E_\ell}| \le 2^{\ell+1} |E_\ell|.
\end{equation}
In view of~\eqref{eq:Fdef} one has 
\bea
F &\subset& \bigcup_{\ell=1}^\infty  \{x\in\R^2\:|\: (\chi_{E_\ell}  \ast \sigma_\delta)(x) \gtrsim \ell^{-2} \,\lambda\} \nn \\
&\subset& \bigcup_{\ell=1}^\infty  \{x\in\R^2\:|\: (\chi_{\wt{E_\ell}} \ast \sigma_{3\delta})(x) \gtrsim 
2^\ell\,\ell^{-2} \,\lambda\}=:\bigcup_{\ell=0}^\infty F_\ell. \label{eq:Fbreak}
\enea
Now fix an arbitrary $\ell\ge1$
and pick a $\delta$-net $\{x_j\}_{j=1}^M\subset F_\ell$. Then
\begin{equation}
\label{eq:Fnet} 
F_\ell \subset \bigcup_{j=1}^M B(x_j,\delta) \subset 
\{x\in\R^2\:|\: (\chi_{\wt{E_\ell}} \ast \sigma_{4\delta})(x) \gtrsim 2^\ell\,\ell^{-2}\lambda\}. 
\end{equation}
Set $\lambda_\ell:=2^\ell\,\ell^{-2}\lambda$.  Since we can assume that $F_\ell\not=\emptyset$, one
concludes from~\eqref{eq:Fnet} that~$\lambda_\ell\gtrsim \delta$. 
By construction, $\wt{E_\ell}$ is discrete on scale $\delta$, i.e., there is a $\delta$-net $\{y_k\}_{k=1}^N\subset \wt{E_\ell}$ with $N\ge1$ so that 
\[ |\wt{E_\ell}| \sim N\delta^2.\]  
By \eqref{eq:Fnet}, every $x_j$ has the property that
\begin{equation}
\label{eq:intersec}
|C^{4\delta}(x_j,1) \cap \wt{E_\ell}| > c_0\, \lambda_\ell\,\delta
\end{equation}
with some absolute constant $c_0$. We will prove that 
\begin{equation}
\label{eq:discr3/2} 
|F_\ell| \lesssim \lambda_\ell^{-3} |\wt{E_\ell}|^2 \text{\ \ or \ \ } M \lesssim \lambda_\ell^{-3}\delta^{2} N^2,
\end{equation}
which implies \eqref{eq:restwt} by summation over $\ell$.
This is a relatively immediate consequence of the fact that two distinct points have
at most two unit circles passing through them. Indeed,  consider the set  
\begin{equation}
\label{eq:Q1def}
 Q= \Bigl\{ (x_j,y_k,y_i)\:\Big|\: \big||x_j-y_k|-1\big|<\delta,\;\big||x_j-y_i|-1\big|<\delta,\; |y_k-y_i|> \frac{c_0}{10}\,\lambda_\ell-\delta \Bigr\}.
\end{equation}
Then
\bea
\card(Q) &\les& N^2 \lambda_\ell^{-1} \label{eq:upper} \\
\card(Q) &\gtrsim& M (\lambda_\ell\, \delta^{-1})^2. \label{eq:lower} 
\enea
The upper bound here comes from the fact that there are at most $\lambda_\ell^{-1}$ many choices of
$\delta$-annuli of radius one passing through two given points at a distance~$\lambda_\ell$ (this is where curvature is used). The lower bound follows from~\eqref{eq:intersec} (here recall that $\lambda_\ell\gtrsim \delta$, or in other words, that there is at least one choice of $y_k$ for every $x_j$). 
Comparison of~\eqref{eq:upper} and~\eqref{eq:lower} yields~\eqref{eq:discr3/2}, 
which in turn yields by summing in~$\ell$, 
\[
 |F|\le \sum_{\ell=1}^\infty |F_\ell| \les \sum_{\ell=1}^\infty \lambda_\ell^{-3} |\wt{E_\ell}|^2 
\les \sum_{\ell=1}^\infty \lambda^{-3}\,2^{-3\ell}\ell^6\,2^{2\ell} |E|^2\les\lambda^{-3} |E|^2. 
\]
This is the desired restricted weak-type form of~\eqref{eq:3/2} (by interpolation it leads to~\eqref{eq:3/2} with an $\epsilon$-derivative loss). 

\noindent This argument generalizes to higher dimensions and surfaces $\{S_x\}_{x\in\Omega}$, with some domain
$\Omega\subset\R^d$, obeying the rotational curvature condition of Phong and Stein~\cite{PhSt}, see also~\cite{Stein1}, which we now recall.

\begin{defi}
\label{def:rotcurv}
Let $\Omega\subset\R^d$ be a region and suppose $\Phi\in C^1(\Omega\times\Omega), \partial^2_{xy} \Phi \in C(\Omega\times\Omega)$. We require that the Monge-Ampere determinant of $\Phi$ satisfies 
\begin{equation}
\label{eq:Steincond}
       \det \left[ \begin{array}{cc} \Phi(x,y) & \partial_y{\Phi}(x,y) \\
                                 \partial_x \Phi(x,y) & \partial^2_{xy} \Phi(x,y) 
               \end{array}
            \right] \not= 0 \text{\ \ whenever\ \ }\Phi(x,y)=0,\;(x,y)\in\Omega^2.
\end{equation}
The hypersurfaces
\[ S_x=\{y\in\Omega\:|\: \Phi(x,y)=0 \} \quad\forall\;x\in\Omega \]
are then said to have nonvanishing rotational curvature. We also impose the following condition
on $\{S_x\}_{x\in\Omega}$: There exists $K=K(\Phi)$ such that 
\begin{equation}
\label{eq:FM} 
 \card\{x\in\Omega \:|\: \Phi(x,y_j)=r_j,\;1\le j\le d\}\les K
\end{equation}
for a.e.~$\{r_j\}_{j=1}^d\in\R^d$. 
\end{defi}

\noindent
Note in particular that~\eqref{eq:Steincond} implies that 
$\min_{x\in\Omega,y\in S_x}|\partial_y \Phi(x,y)|>0$ so that every $S_x$ is 
a submanifold of dimension $d-1$ in~$\Omega$. 
Two examples of such families of surfaces are the unit spheres 
\[ \Phi_S(x,y)=|x-y|-1 \]
and the planes 
\[ \Phi_P(x,y)=x\cdot y-1,\]  see~\cite{Stein1}. The examples of spheres is a special case
of the translation invariant setting, where $S_x=x+S_0$. In that case nonvanishing rotational
curvature is easily seen to be equivalent to nonvanishing Gaussian curvature of~$S_0$. 
For details see Stein~\cite{Stein1} page~494. It is also shown there that 
nonvanishing rotational curvature is a property of $\{S_x\}_{x\in\Omega}$, and does not
depend on the defining function $\Phi$. It is also
 invariant under smooth changes of coordinates in $x$ and~$y$ separately. 
Restricting to compact subsets of~$\Omega$ and changing coordinates we may 
therefore assume that $\Omega=[0,1]^d$ with the understanding that all bounds hold uniformly up to the boundary.
We shall make this assumption for the remainder of this section.
Note that in contrast to previous works we only require the minimal regularity under which~\eqref{eq:Steincond}
makes sense. It is conceivable that Theorem~\ref{thm:rotcurv} below holds in even less regular
situations, but we do not explore this issue here.

We will use the following notation: For any
points $y_j\in S_x$, $1\le j \le d$, the simplex {\em inside the hypersurface} $S_x$ defined by these points 
 is denoted by $\Delta_{S_x}(y_1,\ldots,y_d)$. 
This is well-defined as long as $\{y_j\}_{j=1}^d$ all belong to one
coordinate chart, since one can then define the simplex using coordinates. This of course depends
on the choice of coordinates, but the area of $\Delta_{S_x}(y_1,\ldots,y_d)$ (which is the quantity we are most
interested in) does not change by more than a constant under a change of coordinates. For the case
of spheres, i.e., $\Phi_S(x,y)=|x-y|-1$, our $\Delta_{S_x}(y_1,y_2,y_3)$ is
the spherical triangle spanned by the points $y_1,y_2,y_3$, whereas for the planes given by $\Phi_P(x,y)=x\cdot y-1$ the simplex $\Delta_{S_x}(y_1,y_2,y_3)$ is the Euclidean triangle in the plane defined by $y_1,y_2,y_3$. 
The main goal of this section is to prove the following theorem.

\begin{theorem}
\label{thm:rotcurv}
Let $\{S_x\}_{x\in\Omega}$ be a family of hypersurfaces in $\R^d$, $d\ge2$, as in Definition~\ref{def:rotcurv}.
Define the averaging operators 
\[ Af(x) := \int_{S_x} f(y)\,\sigma_x(dy)\]
where $\sigma_x$ is surface measure on $S_x$. Then one has the restricted weak-type bound
\[ \|Af\|_{L^{d+1,\infty}(\Omega)} \le C \|f\|_{L^{\frac{d+1}{d},1}(\Omega)} \]
where $C$ depends only on $\Phi$ and the dimension $d$.
\end{theorem}

\noindent The proof of the theorem will follow the outline of the previous argument for circles in the plane.
We will require two auxiliary lemmas, the first of which is the higher-dimensional analogue of the fact
that two points have at most two unit circles passing through them (which was used in the previous argument to obtain the upper bound on $\card(Q)$), see the following Lemma~\ref{lem:phong}. 
The second lemma then shows that the condition that we impose in Lemma~\ref{lem:phong} can be made to hold
generically (this is the analogue of the fact that the separation condition on $y_i,y_k$ in~\eqref{eq:Q1def} is harmless). The formulation of the following lemma might appear unnecessarily complicated, due to the presence of the minimum in~\eqref{eq:mesup}. But it is in fact essential for the proof that we define the set in~\eqref{eq:mesup} by means of this minimization, see the remark after the proof.

\begin{lemma}
\label{lem:phong} Let $\Phi$ be as in Definition~\ref{def:rotcurv}.
There exist a small constant $r_0$ and a large constant $C$ depending on $\Phi$ such that for any  
$\{y_j\}_{j=1}^d\in \Omega$ with $\max_{2\le j\le d}|y_j-y_1| < r_0$ one has
\begin{equation}
\label{eq:mesup}
\Bigl|\Bigl\{x\in\Omega\:|\: \max_{1\le i\le d}\dist(S_x,y_i) < \delta,\; 
\min_{\substack{y_i'\in S_x,\,|y_i'-y_i|<C\delta\\1\le i\le d}}   |\Delta_{S_x}(y_1',\ldots,y_d')|>\lambda\Bigr\}\Bigr| \le C\,\delta^d\,\lambda^{-1}
\end{equation}
for all $\delta>0$. 
\end{lemma}
\noindent
\underline{Proof:} Fix $Y:=\{y_j\}_{j=1}^d$ as above and denote the set on the left-hand side 
of~\eqref{eq:mesup} by~$\Lambda(Y)$. Consider the map
\begin{equation}
\label{eq:Vdef} 
V_Y:\Lambda(Y)\to B(0,C_1\,\delta)\subset \R^d,\quad V_Y(x) := \{\Phi(x,y_j)\}_{j=1}^d,
\end{equation}
where $B(0,C_1\delta)\subset \R^d$ is a ball, and $C_1\min_{x,y\in\Omega}|\partial_y\Phi(x,y)| \ge 1$. 
The goal is to estimate $|\Lambda(Y)|$. This will be accomplished by means of the change of variables formula
\[ \int_{\Lambda(Y)} |\det[DV_Y(x)]|\,dx = \int_{\R^d} \card\{x\in\Lambda(Y)\:|\: V_Y(x)=(r_1,\ldots,r_d)\}\,dr_1\,\ldots\,dr_d,\]
see Federer~\cite{Fed}, Theorem~3.2.3. 
In view of~\eqref{eq:Vdef} and~\eqref{eq:FM}, the right-hand side
is 
\[ \les |B(0,C_1\,\delta)|\les \delta^d. \] 
Therefore, we need to prove that 
\begin{equation}
\label{eq:claimDV} |\det[DV_Y(x)]| \gtrsim \lambda.
\end{equation}
on~$\Lambda(Y)$. One has
\bea
\det[DV_Y(x)] &=& \det \left[ \begin{array}{l} \partial_x \Phi(x,y_1) \\
					     \partial_x \Phi(x,y_2) \\
					     ............... \\
					    \partial_x \Phi(x,y_d)
			    \end{array}
                     \right].  \label{eq:DV}
\enea
To prove \eqref{eq:claimDV}, one invokes the following two well-known consequences of the rotational curvature
condition. For all $x\in\Omega$, $y\in S_x$:
\begin{enumerate} 
\item $\partial^2_{xy}\Phi(x,y)$ has maximal rank $d-1$ on the tangent space $T_y(S_x)$. 
\item $\partial_x\Phi(x,y)$ is transverse to the space 
$ W_{xy}:= \span\{\partial^2_{xy}\Phi(x,y)\vv\:|\: \vv\in T_y(S_x)\}.$
\end{enumerate}
These properties follow from the identity
\begin{equation}
\label{eq:rotcurv} 
\left[ \begin{array}{cc} \Phi(x,y) & \partial_y{\Phi}(x,y) \\
                                 \partial_x \Phi(x,y) & \partial^2_{xy} \Phi(x,y) 
               \end{array}
            \right] 
\left[ \begin{array}{c} \alpha \\ \vv
               \end{array}
            \right] 
= 
\left[ \begin{array}{c} 0 \\ \alpha \partial_x \Phi(x,y) + \partial^2_{xy} \Phi(x,y) \vv
               \end{array}
            \right], 
\end{equation}
which holds for all $y\in S_x$ and $\vv\in T_y(S_x)$. Indeed, if i) were to fail, then 
for $\alpha=0$ and some $\vv\in T_y(S_x)$, $\vv\not=0$, the right-hand side of~\eqref{eq:rotcurv} would
vanish, contradicting~\eqref{eq:Steincond}. If~ii) were to fail, then for  
$\alpha=1$ and some $\vv\in T_y(S_x)$ a contradiction would result. 
 i) implies that for $r_0>0$ small (depending only on $\Phi$)  and any $x\in\Omega,\;y_1\in S_x$, the map 
\[ \El_{x,y_1}:\;y\in B(y_1,r_0)\cap S_x \mapsto \partial_x\Phi(x,y) \]
defines a diffeomorphism onto its image, which we denote by $\plates_{x,y_1}$. Moreover, by ii), 
the vector $\partial_x\Phi(x,y_1)$ is transverse to~$\plates_{x,y_1}$ provided $r_0$ is sufficiently small. \\
We now use this fact to show that the absolute value of the determinant in~\eqref{eq:DV} can be
bounded from below by means of a similar determinant in which the $y_j$ are replaced with $y_j'\in S_x$
and $|y_j-y_j'|<C\delta$. Indeed, let $y_j''\in S_x$ with $1\le j\le d$  be such that 
$|y_j-y_j''|=\dist(y_j,S_x)<\delta$ (the latter inequality being part of the definition of~$\Lambda(Y)$). Then
\[ |\partial_x\Phi(x,y_j)-\partial_x\Phi(x,y''_j)|\le \max_{\Omega^2}|\partial_{xy}^2 \Phi(x,y)|\,|y_j-y_j''| \les \delta.\]
It follows from the properties of $\El_{x,y_j''}$ that for a large constant $C$ (depending only on $\Phi$) 
there exist $y_j'\in B(y_j'',C\delta)\cap S_x$ for all $1\le j\le d$ so that 
\[ \partial_x\Phi(x,y_j') \text{\ \ is parallel to\ \ } \partial_x\Phi(x,y_j). \]
Since the lengths of these vectors are comparable, one concludes that
\begin{equation}
\label{eq:DVmod} 
\Bigl|\det[DV_Y(x)]\Bigr| \gtrsim \Bigl| \det \left[ \begin{array}{l} \partial_x \Phi(x,y_1') \\
					     \partial_x \Phi(x,y_2') \\
					     ............... \\
					    \partial_x \Phi(x,y_d')
			    \end{array}
                     \right] \Bigr|, 
\end{equation}
where $y_j'\in S_x$ with $\max_{1\le j\le d}|y_j-y_j'|<C\delta$. This is precisely the situation in
which we can invoke the condition on the volume of the simplex $\Delta_{S_x}(y_1',\ldots,y_d')$ that was included into the definition
of the set~$\Lambda(Y)$. To do so, recall that the vector $\partial_x\Phi(x,y_1')$ is transverse to~$\plates_{x,y_1'}$. 
Therefore, the volume of
the simplex with vertices $0$ and $\partial_x\Phi(x,y_j')$, $1\le j\le d$ (which is the same as the 
determinant in~\eqref{eq:DVmod} up to a constant factor), is comparable to the volume of the simplex
determined by the points $\partial_x\Phi(x,y_j')$, $1\le j\le d$ {\em inside the hypersurface} $\plates_{x,y_1'}$ (see the comments preceding Lemma~\ref{lem:phong}). Finally, under the diffeomorphism~$\El_{x,y_1'}$ this
volume remains comparable to the volume of the simplex $\Delta_{S_x}(y_1',\ldots,y_d')$, which is assumed to be at least~$\lambda$, and we are done.
\eproof

\noindent \underline{Remark:} The reader might wonder if it is necessary to define the set $\Lambda(Y)$ in terms
of the minimum of all $y_i'$ rather than one fixed choice of $y_i'$, say the closest point on~$S_x$ to~$y_i$. 
This is in fact not so, as can be seen from the example of the planes, i.e., with $\Phi_P(x,y)=x\cdot y-1$.
Indeed, take $d=3$ and $x=(0,0,1)$, say. Then it is clear that 
\[ DV_Y(x) =    \left[ \begin{array}{l} y_1 \\
                                          y_2 \\
					  y_3
			    \end{array}
                     \right],
\]
so that $\det DV_Y(x)=0$ if $y_1,y_2,y_3$ are coplanar. It is evident that one can make such a choice of
$y_1,y_2,y_3$ that are $\delta$-close to the horizontal plane $\pi$  at height one, but such that 
the triangle that is obtained by projecting the points onto $\pi$ has nonzero area. On the other hand,
minimizing over all $y_1',y_2',y_3'$ as in the lemma produces a triangle in~$\pi$ with zero area, as desired.
Thus it is necessary to state the condition as in~\eqref{eq:mesup} if one wishes to have the lower bound~\eqref{eq:claimDV} on~$|\det DV_Y|$.

\noindent The following lemma is needed to ensure that the lower bound on the volume of the simplex
required in~\eqref{eq:mesup} holds for a typical choice of $\{y_j\}_{j=1}^d$. The
reader should think of the set $S$ in the following lemma as lying in a single coordinate
patch of the hypersurface~$S_x$ for an arbitrary $x\in\Omega$. This explains the appearance
of~$\R^{d-1}$ in the following lemma. 

\begin{lemma}
\label{lem:randsimp}
Let $S\subset \R^{d-1}$ have positive and finite measure. Let $\Prob=|S|^{-1}\chi_S\,dx$ be
normalized Lebesgue measure restricted to~$S$. Then for all $\eps>0$, 
\begin{equation}
\label{eq:epsprob}
\Prob[(y_1,\ldots,y_d)\in S^d\:|\: |\Delta(y_1,\ldots,y_d)|\le \eps|S|] \le C_d\,\eps
\end{equation}
where $C_d$ is a purely dimensional constant. In particular, 
for a random choice (relative to $\Prob$) of $y_j\in S$ for $1\le j\le d$ the simplex $\triangle(y_{1},y_{2},\ldots,y_{d})$
spanned by these points satisfies
\[ \bigl|\triangle(y_{1},y_{2},\ldots,y_{d})\bigr| > c_0\,|S|\]
with probability at least $\half$. Here $c_0$ is a constant that depends only on the dimension. 
\end{lemma}
\noindent
\underline{Proof:} Clearly,
\bea
&& \Prob[(y_1,\ldots,y_d)\in S^d\:|\: |\Delta(y_1,\ldots,y_d)|\le \eps|S|]  \nn \\
&& = |S|^{-d}\,\int_{S^d} \chi_{\bigl[|(y_2-y_1)\wedge (y_3-y_1) \wedge \ldots \wedge (y_d-y_1)| \le \eps|S|\bigr]}
     \,dy_1\ldots dy_d. \label{eq:probrep}
\enea
Intuitively, one might guess that the case of $S$ being a ball (or equivalently, an ellipsoid since
the events we are considering are affinely invariant)  is the worst. Indeed, if the set is very fragmented,
then typically the volume of a simplex with vertices chosen at random from the set will be much larger. 
This suggests using a rearrangement inequality.
The most general of its kind is the Brascamp, Lieb, Luttinger rearrangement
inequality~\cite{BLL}.  Recall that this inequality states that for linear transformations
$A_j:(\R^n)^k\to\R^n$ and functions $f_j\ge0$ defined on $\R^n$  for $1\le j\le m$,  one has
\begin{equation}
\label{eq:BLL}
   \int\prod_{j=1}^m f_j(A_j(x_1,\ldots,x_k))\,dx_1\ldots dx_k
   \le \int \prod_{j=1}^m f_j^{**}(A_j(x_1,\ldots,x_k))\,dx_1\ldots dx_k,
\end{equation}
where $f_j^{**}$ is the nonincreasing rearrangement of $f_j$, see~\cite{BLL}.
In this form, it does not apply to~\eqref{eq:epsprob} because of the indicator function inside
the integral. However, as observed by Christ~\cite{christ} Theorem~4.2, the proof from~\cite{BLL} carries over
verbatim if an indicator function of a Steiner convex set is inserted on both sides of~\eqref{eq:BLL}.
More precisely, we say that $K\subset (\R^n)^m$ is {\em Steiner convex} (see Definition~4.1 in~\cite{christ})
if for every orthonormal basis $(\nu_1,\nu_2,\ldots, \nu_n)$ of $\R^n$ and every $t\in (\R^{n-1})^m$
the subset
\[ 
\{(x_1,\ldots,x_m)\in K\:|\: (\la x_j,\nu_1\ra,\ldots,\la x_j,\nu_{n-1}\ra) = t_j \text{\ \ for\ \ }1\le j\le m\}
\]
is convex, and balanced in the sense that it is invariant under the mapping
\[ 
(x_1,\ldots,x_m)=\Bigl(\sum_{i=1}^n t^i_1\,\nu_i,\ldots,\sum_{i=1}^n t^i_m\,\nu_i\Bigr)
\to \Bigl(\sum_{i=1}^{n-1} t^i_1\,\nu_i-t_1^n\nu_n,\ldots,\sum_{i=1}^{n-1} t^i_m\,\nu_i-t_m^n\nu_n\Bigr).
\]
The proof in~\cite{BLL} then implies that for such a Steiner convex set~$K$ one has
\begin{equation}
\label{eq:mike}
\int \prod_{j=1}^m f_j(x_j)\,\chi_K(x)\,dx_1\ldots dx_m \le 
\int \prod_{j=1}^m f_j^{**}(x_j)\,\chi_K(x)\,dx_1\ldots dx_m,
\end{equation} 
where $x=(x_1,\ldots,x_m)$ and $f_j\ge0$ are defined on $\R^n$. For the proof of this see 
Theorem~4.2 in~\cite{christ}. It follows from the multilinearity of the determinant
and the invariance of the volume under orthogonal transformations that the set
\[ K_A:=\bigl\{y=(y_1,\ldots,y_d)\in (\R^{d-1})^d\:|\: -A \le (y_2-y_1)\wedge (y_3-y_1) \wedge \ldots \wedge (y_d-y_1) \le A\bigr\} \]
is Steiner convex for every $A$. Thus~\eqref{eq:mike} implies that 
\begin{equation}
\label{eq:prob3} 
|S|^{-d}\,\int \chi_{K_A}(y) \prod_{i=1}^d \chi_S(y_i) \,dy_1\ldots dy_d \le 
|S|^{-d} \, \int \chi_{K_A}(y) \prod_{i=1}^d (\chi_{S})^{**}(y_i) \,dy_1\ldots dy_d
\end{equation}
for every $A$. But $(\chi_{S})^{**}=\chi_{S^{**}}$, where $S^{**}$ is the ball centered at the 
origin with the same volume as~$S$. In particular, setting $A=\eps|S|$ shows that the probability in~\eqref{eq:probrep} is largest for a ball. For the ball it is an easy matter to prove the lemma. 
Firstly, one can take the radius of the ball to be equal to~$1$. In other words,  w.l.o.g.~$S=B(0,1)$. 
Secondly, recall the formula (see Drury~\cite{drury} and~\cite{christ})
\begin{equation}
\label{eq:drury} 
dy_1\ldots dy_{d-1} = c_d |(y_2-y_1)\wedge \ldots (y_{d-1}-y_1)|\,\lambda_\pi(dy_1)\ldots
\lambda_\pi(dy_{d-1})\,d\pi, 
\end{equation}
where $d\pi$ is Haar measure on the hyperplanes $\pi\in M_{d-1,d-2}$ in $\R^{d-1}$, 
and $\lambda_\pi$ is Lebesgue measure on the plane~$\pi$. 
Therefore, the second integral in~\eqref{eq:prob3} is equal to ($\kappa_d$ being another dimensional constant)
\bea
&&  c_d\, \int_{M_{d-1,d-2}} \int_{(B(0,1)\cap\pi)^{d-1}} \int_{B(0,1)} \chi_{\Bigl[\dist(y_d,\pi) < \kappa_d\eps\frac{|B(0,1)|}{|(y_2-y_1)\wedge \ldots (y_{d-1}-y_1)|}\Bigr]}\, dy_d\, \nn \\
&& \qquad\qquad\qquad\,|(y_2-y_1)\wedge \ldots (y_{d-1}-y_1)|\,\lambda_\pi(dy_1)\ldots \lambda_\pi(dy_{d-1})\,d\pi \nn \\
&&  \les  \eps \int_{M_{d-1,d-2}} \int_{(B(0,1)\cap\pi)^{d-1}} \,\lambda_\pi(dy_1)\ldots \lambda_\pi(dy_{d-1})\,d\pi, \nn \\
&& \les \eps. \nn
\enea
as desired. The lemma is proved.
It is natural to ask whether one can give a proof that does not rely on the rearrangement
inequalities. We have worked out such an argument for the case of $d=3$ that is quite short,
but of course it does not show that the ball is the extremal case. The idea is
to work with~\eqref{eq:drury} directly on the set~$S\subset \R^2$, which is completely arbitrary 
(up to having finite and positive measure). Clearly, $M_{2,1}$ is just the space $\lines$ of lines $\ell$ 
in the plane
parameterized by $\ell=(\phi,h)$, where $0\le\phi<\pi$ is the angle the line $\ell$  makes with the horizontal, and with $h$ being the signed distance from~$\ell$ to the origin. Thus~\eqref{eq:drury} now becomes
\[ dy_1dy_2=c_2|t_2-t_1|d\lambda_\ell(t_1)d\lambda_\ell(t_2)\, d\ell = c_2|t_2-t_1|d\lambda_\ell(t_1)d\lambda_\ell(t_2)\, \frac{d\pi}{\pi}\,dh.\]
We shall write $y(\ell,t)$ to denote the point $y\in\ell$ at position $t$ on $\ell$. The choice of origin
$t=0$ on~$\ell$ can be made unique by setting it equal to closest point on~$\ell$ to the origin in~$\R^2$. 
We will use the following notation. With $\proj(w,v)$ being the projection of the vector $v$ onto the unit vector~$w$, define
\[ \strip(y,\phi,a):=\{x\in\R^2\:|\: |\proj(x-y,ie^{i\phi})|< a\}. \]
Thus, $\strip(y,\phi,a)$ is just the strip around a line passing through $y$ with angle $\phi$ and width~$2a$.
We will make use of the following elementary property of triangles $\triangle(y_1,y_2,y_3)$: There is always
a pair of sides, say $y_1 y_2$ and $y_1 y_3$, so that
\[ \overline{y_1 y_2} \le 2 \overline{y_1 y_3} \le 4\overline{y_1 y_2}. \]
Now proceed as follows: 
\bea
&& \Prob(|\triangle(y_1,y_2,y_3)|\le \eps|S|) \le 3\Prob(|\triangle(y_1,y_2,y_3)|\le \eps|S|,\,
\overline{y_1 y_2} \le 2 \overline{y_1 y_3} \le 4\overline{y_1 y_2}) \nn \\
&& = |S|^{-3}\,3c_2\, \int_{\lines} \int_{S\cap\ell}\int_{S\cap\ell} \Bigl|\Bigl\{y_3\in S\:|\: \dist(y_3,\ell)<\frac{2\eps|S|}{|t_1-t_2|},\;|y_3-y_1(t,\ell)|\sim |t_1-t_2|\Bigr\}\Bigr| \nn \\
&& \mbox{\hspace{4.7in}}  |t_1-t_2|\,dt_1dt_2\,d\ell \nn \\
&& \les |S|^{-3}\,\sum_{j\in\Z} 2^{2j} \int_{\lines} \int_{S\cap\ell} \Bigl|\Bigl\{y_3\in S\:|\: \dist(y_3,\ell)<2^{-j+2}\eps|S|,\,|y_3-y_1(t,\ell)|\sim 2^j\Bigr\}\Bigr|\,dt_1d\ell \nn \\
&& \les |S|^{-3}\,\sum_{j\in\Z} 2^{2j} \int_0^\pi \int_S\int_S \chi_{[|y_1-y_3|\sim 2^j]}\, \chi_{\strip(y_3,\phi,2^{-j+2}\eps|S|)}(y_1) \, dy_1dy_3\,d\phi \nn \\
&& \les |S|^{-3}\,\sum_{j\in\Z} 2^{2j}  \int_S\int_S \chi_{[|y_1-y_3|\sim 2^j]}\, \Bigl(1+\frac{|y_1-y_3|}{2^{-j+2}\eps|S|}\Bigr)^{-1}\,dy_3\,dy_1 \nn \\
&& \les |S|^{-3}\,\sum_{j\in\Z} 2^{2j}  \int_S\int_S \chi_{[|y_1-y_3|\sim 2^j]}\, \Bigl(1+\frac{2^{2j}}{\eps|S|}\Bigr)^{-1}\,dy_3\,dy_1 \nn \\
&& \les |S|^{-3}\,\sum_{j\in\Z} \sup_{k\in\Z}2^{2k}\Bigl(1+\frac{2^{2k}}{\eps|S|}\Bigr)^{-1}  \int_S\int_S \chi_{[|y_1-y_3|\sim 2^j]}\, \,dy_3\,dy_1 \nn \\
&& \les |S|^{-3}\,\eps|S| \int_S\int_S \,dy_3\,dy_1 = \eps, \nn
\enea
as desired.\eproof 

\noindent The previous lemma is insufficient for our purposes because of the minimization
condition that appears in~\eqref{eq:mesup}. However, the following corollary to 
Lemma~\ref{lem:randsimp} addresses this issue. It is formulated for sets which are discrete on scale~$\delta$, 
which is precisely the situation that arises in the proof of Theorem~\ref{thm:rotcurv}.

\begin{cor}
\label{cor:bigsimp}
Suppose
\begin{equation}
\label{eq:Ssep} 
S=\bigcup_{j=1}^N B(x_j,3\delta)\subset \R^{d-1} \text{\ \ where\ \ }
|x_i-x_j|>6\delta \qquad\forall\; 1\le i\not=j \le N.
\end{equation}
Then with probability at least $\frac{6}{(12)^{d}}$, a simplex $\Delta(y_1,\ldots,y_d)$ with
vertices $\{y_j\}_{j=1}^d$ chosen at random from the set~$S$ has the property that
\begin{equation}
\label{eq:minstuff}
 \min_{|y_i-\wt{y_i}|<\frac{\delta}{4},\;1\le i\le d} |\Delta(\wt{y_1},\ldots,\wt{y_d})|  > c_0\,3^{-d+1}|S|
\end{equation}
where $c_0$ is the constant from Lemma~\ref{lem:randsimp}.
\end{cor}
\noindent
\underline{Proof:} For the purposes of this proof only we define a {\em cone}
to be any rotation and translation of the set $\R_+^{d-1}$. 
Given an arbitrary simplex $\Delta(p_1,\ldots,p_d)\subset\R^{d-1}$ it is quite evident that for every 
$1\le j\le d$ there exists a cone $\Gamma_j\subset\R^{d-1}$ with vertex at~$p_j$  such that 
\[ \inf_{p_j'\in \Gamma_j,\,1\le j\le d}|\Delta(p_1',\ldots,p_d')| \ge |\Delta(p_1,\ldots,p_d)|.\]
We will also need the following elementary geometry fact: if $p\in B(0,\delta)$ then 
for any cone with vertex at~$p$ one has
\begin{equation}
\label{eq:goodball} \Gamma \cap B(0,3\delta)\setminus B(0,\delta) \supset B(p',\delta/2) 
\end{equation}
for some point $p'$. The lemma follows easily from these two properties. Indeed,
define $S'=\bigcup_{j=1}^N B(x_j,\delta)$ and apply Lemma~\ref{lem:randsimp}
to the set~$S'$. Then with probability at least~$\half$, a simplex $\Delta(y_1',\ldots,y_d')$ with
vertices chosen at random from the set~$S'$ has the property that
\begin{equation}
\label{eq:start} 
|\Delta(y_1',\ldots,y_d')|  > c_0|S'|= c_0\,3^{-(d-1)}|S|,
\end{equation}
which is the lower bound in~\eqref{eq:minstuff}. Now fix $y_1',\ldots,y_d'\in S'$ as in~\eqref{eq:start} and 
associate with each $y_j'$ its ball $B_j\subset B(x_{k(j)},3\delta)\setminus B(x_{k(j)},\delta)$ 
of radius $\delta/2$
as in~\eqref{eq:goodball}. Here $k(j)$ is determined by $y_j'\in B(x_{k(j)},\delta)$. Let $B_j'\subset B_j$ 
be the ball with the same center as $B_j$ but radius $\delta/4$. By construction, 
for every choice of $y_j\in B_j'$ one has 
\begin{equation}
\label{eq:select}
  \min_{|y_i-\wt{y_i}|<\frac{\delta}{4},\;1\le i\le d} |\Delta(\wt{y_1},\ldots,\wt{y_d})| 
  \ge |\Delta(y_1',\ldots,y_d')|  > c_0\,3^{-d+1}|S|. 
\end{equation}
Observe that no $B_j'$ can come from two different balls $B(x_{k_1},\delta)$, $B(x_{k_2},\delta)$
because of $| x_{k_1} - x_{k_2}|>6\delta$. 
Thus a moment's reflection shows that the passage from $\{y_j'\}_{j=1}^d\subset S'$ to $y_j\in B_j'$ 
decreases the probability (which was $\half$) by at most a factor of
\[ \frac{|B(0,\delta/4)|}{|B(0,\delta)|} \frac{|B(0,\delta)|}{|B(0,3\delta)|} = (12)^{-(d-1)}\]
so that the total probability of picking $y_j\in S$ for $1\le j\le d$ for which~\eqref{eq:select} holds
is at least $\half (12)^{-(d-1)}=6\, (12)^{-d}$. 
\eproof

\noindent
\underline{Remark:} If $S=\bigcup_{j=1}^N B(x_j,3\delta)$, but  the separation
condition in~\eqref{eq:Ssep} does not hold, then one can still apply the lemma at the cost of  reducing the probability
and the lower bound in~\eqref{eq:minstuff} by another dimensional constant. Indeed, simply pass to 
a subset $\wt{S}\subset S$ that does satisfy~\eqref{eq:Ssep} and so that $|\wt{S}|$ is comparable to~$|S|$.

\medskip\noindent
\underline{Proof of Theorem~\ref{thm:rotcurv}:} We need to show that there exists a constant~$C$ 
only depending on $\Phi$ and the dimension~$d$ such that for any $E\subset \Omega=[0,1]^d$ one has
\[ |\{x\in\Omega\:|\: (A\chi_E)(x)>\lambda\}| \le C\, \lambda^{-(d+1)}\,|E|^{d}\]
for all $\lambda>0$. This is equivalent to showing that
\begin{equation}
\label{eq:wtd} 
|\{x\in\Omega\:|\: (A^\delta\chi_E)(x)>\lambda\}| \le C\, \lambda^{-(d+1)}\,|E|^{d}
\end{equation}
uniformly in $\delta>0$ where 
\[ (A^\delta f)(x) = \int_\Omega f(y)\, d\sigma_x^\delta(y),\]
$\sigma^\delta$ being normalized measure on a $\delta$-neighborhood of the hypersurface~$S_x$, which we denote
by~$S_x^\delta$. 
Fix some small $\delta>0$, as well as some $E\subset\Omega $ and $\lambda>0$. Define
\begin{equation}
\label{eq:Fdefd} 
F:=\{x\in \Omega\:|\: (A^\delta\chi_E)(x)>\lambda\}.
\end{equation}
As in the argument dealing with circles, we discretize on scale $\delta$. More precisely,
partition $\Omega$ into squares $\{Q_j\}$ of side-length $\delta$ and let
\[ E_\ell = \bigcup_{j\::\:2^{-\ell}\delta^d < |Q_j\cap E| \le 2^{-\ell+1}\delta^d} Q_j\cap E
 \text{\ \ \ and\ \ \ }\wt{E_\ell} = \bigcup_{j\::\:2^{-\ell}\delta^d<|Q_j\cap E|\le 2^{-\ell+1}\delta^d} Q_j\]
for $\ell\ge1$. Clearly, $E=\bigcup_\ell E_\ell$ and $|\wt{E_\ell}| \sim 2^{\ell} |E_\ell|$. 
By \eqref{eq:Fdefd} one has 
\bea
F &\subset& \bigcup_{\ell=1}^\infty  \{x\in\Omega\:|\: \big(A^\delta\,\chi_{E_\ell} \big)(x) \gtrsim \ell^{-2} \,\lambda\} \nn \\
&\subset& \bigcup_{\ell=1}^\infty  \{x\in\Omega\:|\: \big( A^{C_3\delta}\chi_{\wt{E_\ell}} \big)(x) \gtrsim 
2^\ell\,\ell^{-2} \,\lambda\}=:\bigcup_{\ell=0}^\infty F_\ell, \label{eq:Fteil}
\enea
where $C_3$ is some constant depending on the dimension and $\Phi$ (one can take $C_3\sim 1+\sqrt{d}$).
Now fix an arbitrary $\ell\ge1$ and pick a $\delta$-net $\{x_j\}_{j=1}^M\subset F_\ell$. Then
\begin{equation}
\label{eq:Fnetd} 
F_\ell \subset \bigcup_{j=1}^M B(x_j,\delta) \subset 
\{x\in\Omega\:|\: \big(A^{C_4\delta}\chi_{\wt{E_\ell}} \big)(x) \gtrsim 2^\ell\,\ell^{-2}\lambda\} 
\end{equation}
where $C_4$ is a constant that depends only on $\Phi$ and~$d$. 
Set $\lambda_\ell:=2^\ell\,\ell^{-2}\lambda$.  Since we can assume that $F_\ell\not=\emptyset$, one
concludes from~\eqref{eq:Fnetd} that~$\lambda_\ell\gtrsim \delta^{d-1}$. 
By construction, $\wt{E_\ell}$ is discrete at scale $\delta$, i.e., there is a 
$\delta$-net $\{y_k\}_{k=1}^N\subset \wt{E_\ell}$ with $N\ge1$ so that $|\wt{E_\ell}| \sim N\delta^d$. 
By~\eqref{eq:Fnetd}, every $x_j$ has the property that
\begin{equation}
\label{eq:intersecd}
|S^{C_4\delta}_{x_j} \cap \wt{E_\ell}| > c_5\, \lambda_\ell\,\delta
\end{equation}
with some (small) constant $c_5$. We will prove that 
\begin{equation}
\label{eq:discr3/2d} 
|F_\ell| \lesssim \lambda_\ell^{-d-1} |\wt{E_\ell}|^d \text{\ \ or \ \ } M \lesssim \lambda_\ell^{-d-1}
\delta^{d(d-1)} N^d,
\end{equation}
which implies \eqref{eq:wtd} by summation over $\ell$. 
To prove~\eqref{eq:discr3/2d} we need to apply Corollary~\ref{cor:bigsimp}. In order to do so it will be convenient to ``fatten up'' the set $\wt{E_\ell}$ as follows.
With every $y\in \wt{E_\ell}$ include the entire ball $B(y,C\delta)$ in~$\wt{E_\ell}$ where $C$ is the constant
from Lemma~\ref{lem:phong}. This means, of course, that with every point $y_k$ in the net of $\wt{E_\ell}$
we include all its $C\delta$-neighbors into the net as well. Clearly, this only has the effect of loosing another constant in~\eqref{eq:discr3/2d}, but otherwise everything remains unchanged. 
With this in mind consider the set  
\bea
\label{eq:Q1defd}
 Q &=& \Bigl\{ (x_j,y_{k_1},\ldots,y_{k_d})\:\Big|\: \dist(S_{x_j},y_{k_i})<\delta,\;1\le i\le d,\; \nn \\
&& \mbox{\hspace{2in}} \min_{\substack{y_{k_i}'\in S_{x_j},\,|y_{k_i}'-y_{k_i}|<C\delta\\1\le i\le d}} 
|\Delta_{S_x}(y_{k_1}',\ldots,y_{k_d}')| > c_6\,\lambda_\ell  \Bigr\} \nn
\enea
where $c_6$ is a small constant that depends on $c_5$ and the constants from Corollary~\ref{cor:bigsimp}
(see also the remark following Corollary~\ref{cor:bigsimp}). Then
\bea
\card(Q) &\les& N^d \lambda_\ell^{-1} \label{eq:upperd} \\
\card(Q) &\gtrsim& M (\lambda_\ell\, \delta^{-(d-1)})^{d}. \label{eq:lowerd} 
\enea
The upper bound is an immediate consequence of Lemma~\ref{lem:phong}. 
The lower bound follows from~\eqref{eq:intersecd}, Corollary~\ref{cor:bigsimp} and the remark following it
(apply those with $C\delta$ instead of $\delta$).  
Comparison of~\eqref{eq:upperd} and~\eqref{eq:lowerd} yields~\eqref{eq:discr3/2d}, 
which in turn implies by summation in~$\ell$, 
\[
 |F|\le \sum_{\ell=1}^\infty |F_\ell| \les \sum_{\ell=1}^\infty \lambda_\ell^{-d-1} |\wt{E_\ell}|^d 
\les \sum_{\ell=1}^\infty \lambda^{-d-1}\,2^{-(d+1)\ell}\ell^{2(d+1)}\,2^{d\ell} |E|^d \les \lambda^{-d-1}\,|E|^d, 
\]
and the theorem is proved.\eproof

\section{A geometric proof of a Strichartz estimate}
\label{sec:strich}

Let $\sigma_\delta$ be the normalized measure on the $\delta$ neighborhood of the cone 
$\Gamma\cap \{1<|t|<2\}$. For given $\lambda>0$ and $E\subset [0,1]^2$, let
\[ F=\{(x,t)\in \R^2_x\times \R_t\:|\: (\sigma_\delta \ast \chi_E) (x,t) >\lambda \}, \]
where $\ast$ stands for convolution in $\R^2_x$. We will show below that for any $\eta>0$ there exists a 
constant~$C_\eta$ so that
\begin{equation}
\label{eq:level}
\lambda |F|^{\frac16} \le C_\eta\, \delta^{-\eta}\, |E|^{\frac12}.
\end{equation}
In other words, one has the restricted weak-type bound
\[ \|\sigma_\delta\ast f\|_{L^{6,\infty}(\R^2\times[1,2])} \le C_\eta\,\delta^{-\eta}\,\|f\|_{L^{2,1}(\R^2)}.\]
As usual, this implies  the strong-type bound
\[
\|\sigma_\delta\ast f\|_{L^{6}(\R^2\times[1,2])} \le C_\eta\,\delta^{-\eta}\,\|f\|_{L^{2}(\R^2)}
\]
by means of interpolation with an easy $2\to 6$ bound with a loss~$\delta^{-1}$, say. 
It is now clear that one also has the bound
\begin{equation}
\label{eq:2to6} 
\|\sigma_\Gamma\ast f\|_{L^{6}(\R^2\times[1,2])} \le C_\eps\,\|f\|_{W^{2,\eps}(\R^2)}
\end{equation}
where $\sigma_\Gamma$ is the surface measure on the cone segment $\Gamma=\{|\xi|=t\:|\: 1\le t\le 2\}$. 
For the sake of completeness, we show that~\eqref{eq:2to6} implies a Strichartz-type bound of the form
\begin{equation}
\label{eq:strichloct}
\|u\|_{L^6(\R^2\times [1,2])} \le C_{\epsilon}\, \|f\|_{W^{2,\half+\eps}(\R^2)}.
\end{equation}
for solutions of the wave equation 
\[ \Box u=0,\qquad u|_{t=0}=f,\; \partial_t u|_{t=0}=0. \]
Arguments of this type appear in \cite{KW} and~\cite{W}, but here we proceed somewhat differently. 
Firstly, one writes
\[ (\sigma_\Gamma\ast f)(\cdot,t) = (m_t\hat{f})^{\vee}(\cdot)\chi_{[1,2]}(t)\]
with 
\[ m_t(\xi) = \omega_{+}(t|\xi|)e^{it|\xi|} + \omega_{-}(t|\xi|)e^{-it|\xi|} .\]
This latter representation comes of course from the Fourier transform of the surface measure of the circle.
Thus, 
\begin{equation}
\label{eq:mpm} 
\bigl|\frac{d^k}{dr^k}\,\omega_{\pm}(r)\bigr| \le C_k\,(1+r)^{-\half-k},\;k\ge0,\quad
\lim_{r\to\infty} \sqrt{r}\,\omega_{\pm}(r)=c_{\pm}\not=0.
\end{equation}
Let $R_1, R_2$ be the usual Riesz transforms on $\R^2$ with multipliers $-i\frac{\xi_1}{|\xi|}$, and
$-i\frac{\xi_1}{|\xi|}$, respectively. Then one has the identity (with $\FT$ denoting the Fourier transform)
\bea
&& \FT R_1[(-ix_1\sigma_\Gamma)\ast f]+ \FT R_2[(-ix_2\sigma_\Gamma)\ast f] = t[\omega_{+}(t|\xi|)e^{it|\xi|}-
\omega_{-}(t|\xi|)e^{-it|\xi|}]\chi_{[1,2]}(t)\,\hat{f}(\xi) \nn \\
&& \qquad -it[\omega_{+}'(t|\xi|)e^{it|\xi|}+\omega_{-}'(t|\xi|)e^{-it|\xi|}]\chi_{[1,2]}(t)\,\hat{f}(\xi). \nn
\enea
Hence
\begin{equation}
\label{eq:repTR} 
\int e^{i[x\cdot\xi+t|\xi|]}\omega_{+}(t|\xi|) \hat{f}(\xi)\,d\xi \,\chi_{[1,2]}(t) = 
\half\Bigl((\sigma_\Gamma\ast f)(\cdot,t)+T_tf\Bigr) + E_t f, 
\end{equation}
where 
\bea 
T_t f &=& \frac{1}{t}\,R_1[(-ix_1\sigma_\Gamma)\ast f]+\frac{1}{t}\,R_2[(-ix_2\sigma_\Gamma)\ast f] \nn \\
\widehat{E_t f}(\xi) &=& \half i [\omega_{+}'(t|\xi|)e^{it|\xi|}+\omega_{-}'(t|\xi|)e^{-it|\xi|}]\chi_{[1,2]}(t).\nn
\enea
Then, using this representation and the Sobolev embedding theorem to control the ``error term'' $E_t$, 
\bea
&& \Bigl\| \int e^{i(x\cdot\xi+t|\xi|)} \omega_{+}(t|\xi|)\chi_{[1,2]}(t)\,\hat{f}(\xi)\,d\xi \Bigr\|_{L^6_{x,t}} \nn \\
&& \les \|\sigma_{\Gamma}\ast f\|_{L^6_{x,t}}+ \|(x_1\sigma_{\Gamma})\ast f\|_{L^6_{x,t}} + \|(x_2\sigma_{\Gamma})\ast f\|_{L^6_{x,t}} + \sup_{1\le t\le2} \|(-\triangle)^{\frac13}(E_tf)\|_{L^2}, \nn
\enea
and similarly for $e^{i[x\cdot\xi-t|\xi|]}$. It is easy to see from \eqref{eq:2to6} that for functions $f$ with $\supp(f)\subset [0,1]^2$ 
all the terms on the right-hand side are no larger than~$\|f\|_{W^{2,\eps}(\R^2)}$. Invoking the Mikhlin theorem
to remove $\omega_{+}$ yields
\[ \Bigl\| \int_{\R^2} e^{ix\cdot\xi} \cos(t|\xi|) \hat{f}(\xi)\,d\xi \Bigr\|_{L^6_{x,t}(\R^2\times[1,2])} 
\le C_\eps\, \|f\|_{W^{2,\half+\eps}(\R^2)}\]
provided $\supp(f)\subset [0,1]^2$. This latter condition is now eliminated by means of the finite propagation speed for the wave equation. 

\noindent We will give two different purely geometric-combinatorial proofs of~\eqref{eq:level}, 
see Corollary~\ref{cor:main} and Lemma~\ref{lem:three} below. The setup is the same as in~\cite{tomrev}.
Let $C(x,r)$ be as in the previous section with $x\in[0,1]^2$, and define
\[ \Delta(C(x,r),C(y,s))=\bigl||x-y|-|r-s|\bigr|,\quad d(C(x,r),C(y,s))=|x-y|+|r-s|.\]
 $\Delta$ measures the extent to which $C(x,r),C(y,s)$ are internally tangent. As always, we will need 
to know the area of intersection of two annuli, as well as the diameter of their intersection. 
Writing $\Delta$ and $d$ without arguments for simplicity, 
\bea
\bigl|C^\delta(x,r)\cap C^\delta(y,s)\bigr| &\les& \frac{\delta^2}{\sqrt{(\Delta+\delta)(d+\delta)}} \label{eq:durch} \\
\diam \bigl(C^\delta(x,r)\cap C^\delta(y,s)\bigr) &\les& \sqrt{\frac{\Delta+\delta}{d+\delta}}, \nn
\enea
see Lemma~3.1 in~\cite{tomrev}, for example. A central role is played by the multiplicity function
\[ \mu_\delta^{\C}=\sum_{C\in\C} \chi_{C^\delta}. \]
Given a family $\C$ of circles, we use the notation
\begin{equation}
\label{eq:CCetdef} 
\CCet = \bigl\{\Cb\in\C\:|\: \eps-\delta\le\Delta(C,\Cb)\le 2\eps,\; t/2\le d(C,\Cb)\le  t,\; C^\delta\cap \Cb^\delta\not=\emptyset \bigr\}.
\end{equation}
Thus, $\CCet$ is the collection of circles that are $\eps$-tangent to $C$ and have distance $t$ from~$C$. 
We will always assume that $\C$ is $\delta$-separated in the sense that $d(C,C')>\delta$ for any distinct
$C,C'\in \C$. 
Notice that if $\Cb\in \CCet$, then $|x-\bar{x}|\sim t$. The following lemma shows how to 
reduce~\eqref{eq:level} to a bound on the multiplicity function, cf.~\cite{S1}, \cite{S2} 
for similar statements.

\begin{lemma} 
\label{lem:equiv}
The following are equivalent:
\begin{enumerate}
\item Let $\eta>0$ be arbitrary. Then for $\delta>0$ sufficiently small depending on $\eta$, the estimate~\eqref{eq:level} holds for all $E\subset [0,1]^2$ and $0<\lambda<1$. 
\item Given a $\delta$-separated three parameter family of circles $\C$, one has: For any $\eta>0$ and $\delta$ sufficiently small depending on~$\eta$, there exists $\A\subset\C$ such that $|\A|>\half|\C|$ and
\begin{equation}
\label{eq:vielfach} 
|\{C^\delta\:|\: \mu_\delta^{\A} > \delta^{-\eta}\lambda^{-1}|\A|^{\frac23}\}| \le \lambda|C^\delta|
\end{equation}
for all $C\in\A$, $0<\lambda<1$.
\item Same statement as in ii) but with $|\A|>c_\eta\,|\C|$ for some small constant $c_\eta$. 
\end{enumerate}
\end{lemma}
\noindent
\underline{Proof:} We first deal with iii) implies i). Let $\{(x_j,t_j)\}_{j=1}^N$ be a maximally $\delta$-separated
set of points inside~$F$. It is easy to see that $|F|\sim N\delta^3$. Denote the family $\{C(x_j,t_j)\}_{j=1}^N$ of circles by~$\C$. Then for all $C\in\C$  
\[ |C^\delta\cap E| > \lambda|C^\delta|\]
by definition. Take a small $\eta$ and pass to $\A\subset\C$ as in Corollary~\ref{cor:main}. Fix any $0<\lambda<1$. 
On the one hand,
\[ 
\int_{E:\mu^{\A}_\delta<\delta^{-\eta}\lambda^{-1}|\A|^{\frac23}} \mu^{\A}_\delta(x) \,dx \le \delta^{-\eta}\lambda^{-1}N^{\frac23}|E|.
\]
On the other hand, by \eqref{eq:mult} with $\frac{\lambda}{2}$, 
\[
\int_{E:\mu^{\A}_\delta<\delta^{-\eta}\lambda^{-1}|\A|^{\frac23}} \mu^{\A}_\delta(x) \,dx \ge |\A|\frac{\lambda}{2}\delta.
\]
Thus,
\[ \delta^{-\eta}|E|\gtrsim \lambda^2 N^{\frac13}\delta \sim \lambda^2 |F|^{\frac13},\]
as desired.

\noindent To show that i) implies ii) we use an argument from~\cite{S1},\cite{S2}.
First note that the dual statement to (the strong form of) \eqref{eq:level} is
\begin{equation}
\label{eq:dual}
\|\sum_{j=1}^N a_j \chi_{C^\delta(x_j,t_j)} \|_{L^2(\R^2)} \les \delta^{-\eta}\, \sqrt{\delta} 
(\sum_{j=1}^N |a_j|^{\frac65})^{\frac56}
\end{equation}
for any $\delta$-separated three parameter family of circles. 
Now suppose that for at least half the circles in a given family~$\C$ one has
\[ |\{C^\delta\:|\: \mu^{\C} > \delta^{-\eta} \lambda^{-1} |\C|^{\frac23} \}| > \lambda|C^\delta| \]
for some choice of $0<\lambda<1$. 
Pigeonholing as usual we may assume that $\lambda$ is fixed for all circles $C\in\B$ with $|\B|\gtrsim |\log\delta|^{-1}|\C|$. 
Then set
\[ E=\{ \mu^{\C} > \delta^{-\eta} \lambda^{-1} |\C|^{\frac23} \} .\]
Applying both our assumption i) (in strong form, say) and the dual~\eqref{eq:dual} with this choice of~$E$
and $N=|\C|$ yields 
\bea
 \lambda (N\delta^3)^{\frac16} &\les& \delta^{-\frac{\eta}{10}}|E|^{\frac12} \nonumber \\
 \delta^{-\eta} \lambda^{-1} N^{\frac23} |E|^{\frac12} &\les& \delta^{-\frac{\eta}{10}}\sqrt{\delta} N^{\frac56}, \nonumber
\enea
which are incompatible.
\eproof

\noindent Since the Strichartz estimate \eqref{eq:26} assures that i) holds, this lemma proves that~ii), iii) also hold. It is common knowledge that the Strichartz estimate under consideration is (only) optimal for the so called Knapp example, i.e., a slab of dimensions $1\times\delta\times\sqrt{\delta}$ that lies on a light cone.
In our setting this would correspond to a family of circles $\C$ with $|\C|\sim \delta^{-\frac32}$. Note that
in this case the bound in~\eqref{eq:vielfach} is optimal with $\lambda=\sqrt{\delta}$, as expected.\\
In what follows, we give two different direct proofs of ii) and~iii). 
The first one is based on Marstrand's three circle lemma~\cite{Mar}. This is a continuum analogue of the circles of Appolonius. We will not repeat the heuristics for these ideas, as they can be found in~\cite{tomrev}.  For the convenience of the reader we do however reproduce the statement of the three circle lemma from~\cite{tomrev}.

\begin{lemma} 
\label{lem:mar}
With some sufficiently large numerical constant $a_0$, assume that $\epsilon,t,\lambda\in
 (0,1)$ satisfy $a_0\, \epsilon \leq t\lambda^2$.  Fix three 
circles $C(x_i,r_i), 1\leq i\leq 3$. Then for $\delta\leq\eps$
 the set
\bea
\overline{\Omega}_{\epsilon t\lambda} &:=& \{(x,r)\in\R^2\times\R \:|\:
\Delta(C(x,r),C(x_i,r_i))<\epsilon\quad\forall i, \nn \\
&& d(C(x,r),C(x_i,r_i))>t\;\;\forall i, \quad C_{\delta}(x,r)\cap 
C_\delta(x_i,r_i)\neq\emptyset\;\;\forall i, \nn\\ 
&& \dist(C_{\delta}(x,r)\cap
 C_\delta(x_i,r_i),C_{\delta}(x,r)\cap C_\delta(x_j,r_j))\geq\lambda\;
\;\;\forall i,j:i\neq j\} \nn
\enea
is contained in the union of two ellipsoids in $\R^3$ each
 of diameter $\lesssim\frac{\eps}{\lambda^2}$ and volume 
$\lesssim\frac{\eps^3}{\lambda^3}$.
\end{lemma}

\noindent We now show how this lemma immediately leads to the desired multiplicity bound. The case distinction
that arises in the proof has to do with the degenerate configuration where three circles are tangent at one
point. 

\begin{lemma}
\label{lem:three}
Let $\C$ be a $\delta$-separated three parameter family of circles. If $\delta$ is sufficiently small, then 
there exists $\A\subset\C$ such that $|\A|>\half|\C|$ and 
\begin{equation}
\label{eq:mult2} 
|\{C^\delta\:|\: \mu_\delta^{\C} > |\log\delta|^{5}\lambda^{-1}|\C|^{\frac23}\}| \le \lambda|C^\delta|
\end{equation}
for all $C\in\A$, $0<\lambda<1$.
\end{lemma}
\noindent
\underline{Proof:} Suppose for at least half the circles in $\C$ one has
\[
|\{C^\delta\:|\: \mu_\delta^{\C} > |\log\delta|^{5}\lambda^{-1}|\C|^{\frac23}\}| > \lambda|C^\delta|
\]
for some choice of $\lambda$. Pigeonholing yields $\delta\le\eps\les t\le 1$,  $0<\lambda<1$, and $\B\subset\C$ with $|\B|\gtrsim |\log\delta|^{-3}|\C|$  such that 
\begin{equation}
\label{eq:high}
|\{C^\delta\:|\: \mu_\delta^{\CCet} > |\log\delta|^{2}\lambda^{-1}|\C|^{\frac23}\}| > \lambda|C^\delta|
\end{equation}
for all $C\in\B$. We now distinguish two cases. For convenience we denote the set on the left--hand side of~\eqref{eq:high} by~$H(C^\delta)$ (the ``high multiplicity part'' of~$C^\delta$). 

\noindent \underline{Case 1:} For all $C\in\B$ and all $x\in\R^2$ one has
\begin{equation}
\label{eq:twoends}
 |H(C^\delta)\cap B(x,a_0\sqrt{\eps/t})| \le  \frac{\lambda}{100}|C^\delta|,
\end{equation}
$a_0$ being the constant from Lemma~\ref{lem:mar}. 
Let $R\ge a_0\sqrt{\eps/t}$ be maximal with the property that \eqref{eq:twoends} holds
with $R$ instead of~$\sqrt{\eps/t}$, i.e., for all $C\in\B$ and all $x\in\R^2$ one has
\begin{equation}
\label{eq:twoendsmax} |H(C^\delta)\cap B(x,R)| \le  \frac{\lambda}{100}|C^\delta|.
\end{equation}
Consider the set
\bea 
 Q &:=& \{(C,C_{i_1},C_{i_2},C_{i_3})\:|\: C\in\B,\; C_{i_1},C_{i_2},C_{i_3}\in\CCet,\nn \\ 
 && \qquad \dist(C\cap C_{i_\ell}, C\cap C_{i_k}) \ge R \text{\ \ for all\ }1\le k<\ell\le 3. \}. \nn
\enea
In that case one can apply the three circle lemma to conclude that
\[ 
|\B|\bigl(|\log\delta|^{2}\lambda^{-1}|\C|^{\frac23}\frac{\lambda\sqrt{\eps t}}{\delta}\bigr)^3 \les \frac{\eps^3}{\delta^3}R^{-3} |\C|^3.
\]
The upper bound follows from Lemma~\ref{lem:mar}, whereas the lower bound follows from~\eqref{eq:twoendsmax}, \eqref{eq:high} and~\eqref{eq:durch} (apply the latter with $\Delta=\eps, d=t$ and conclude that the number of 
curvilinear rectangles of area $\delta^2/\sqrt{\eps t}$ that are each hit by about $\mu=|\log\delta|^{2}\lambda^{-1}|\C|^{\frac23}$ annuli is~$\frac{\lambda\delta}{\delta^2/\sqrt{\eps t}}$). Simplifying the previous inequality yields
\[ |\log\delta|^3 (\eps t)^{\frac32} \les \eps^3 R^{-3},\]
which contradicts $R\ge a_0\sqrt{\eps/t}$.

\noindent \underline{Case 2:} For one $C_0\in\B$ there is an $x_0\in\R^2$ for which \eqref{eq:twoends} fails.
In that case, we simply compare the number of circles that actually do intersect $C_0^\delta$ inside the ball
$B(x_0,a_0\sqrt{\eps/t})$ because of our multiplicity assumption to the largest possible number that can 
intersect it there. 
With $\mu=|\log\delta|^{2}\lambda^{-1}|\C|^{\frac23}$ this yields 
\begin{equation}
\label{eq:arc}
 \mu\lambda\frac{\sqrt{\eps t}}{\delta} \les \min(|\C|,\sqrt{t\eps}\;t\eps\,\delta^{-3}) \les |\C|^{\frac23} \frac{\sqrt{\eps t}}{\delta},
\end{equation}
which implies $\mu\les \lambda^{-1}|\C|^{\frac23}$, a contradiction.
The right--hand side of \eqref{eq:arc} comes from the fact that the centers of the circles contributing
to~$\mu$ belong to a rectangle of dimension $t\times \sqrt{\eps t}$, the freedom in the radius then
giving another~$\eps$.
\eproof

\noindent We would like to emphasize that this argument carries over verbatim to the case of averages
over $\delta$-neighborhoods of curves satisfying Sogge's cinematic curvature condition. This is due
to the fact that Lemma~\ref{lem:mar} was shown to hold in this context by Kolasa and Wolff~\cite{KW}.
Consequently, the proof of Lemma~\ref{lem:three}  yields the estimate (cf.~\eqref{eq:2to6})
\[ \|{\mathcal A} f\|_{L^{6}(\R^2\times[1,2])} \le C_\eps\,\|f\|_{W^{2,\eps}(\R^2)} \]
where
\[ {\mathcal A}f(x,r)=\int f(y)\,d\sigma_{\gamma_{x,r}}(y) \]
and $\gamma_{x,r}$ is a family of curves with cinematic curvature.\\
We now present a different proof of \eqref{eq:level} that does not rely on the three circle lemma.
Rather, it relies on a ``two circle'' lemma from~\cite{S2}. This is  a device to
control the number of $\delta$-separated  circles  that are tangent to two given ones.
In contrast to the situation of Lemma~\ref{lem:mar}  all circles that are tangent to two given ones form a one
parameter family. Note that in a purely combinatorial setting it is meaningless to work with such a 
two circle lemma, as all circles could belong to this family. It turns out, however, that in our context in
which $\delta$-separateness is imposed, such a device turns out to be useful.  
Hence this is an example of a method that works only for continuum incidence problems, but has no
analogue in incidence geometry per se. We start be recalling this two-circle lemma from~\cite{S2}.

\begin{lemma} 
\label{lem:twocirc}
Suppose  $C_2\in\C^{C_1}_{\beta\tau}$. Then 
\begin{equation}
\label{eq:int}
|\durch|\les\frac{\eps t^2}{\delta^3}\min(\sqrt{\frac{\epsilon}{\tau}},\frac{\epsilon}{\sqrt{\beta \tau}}).
\end{equation}
\end{lemma}
\noindent
\underline{Proof:} As the details are exactly the same as in~\cite{S2}, we do not repeat them. 
The bound~\eqref{eq:int} is actually proved there implicitly, see Lemma~2.5 in that paper. 
The main difference is that~\cite{S2} works with a two-parameter family of circles in the plane, whereas here we need to consider a three parameter family. This, however, only requires changing various bounds in~\cite{S2}
by a factor of $\frac{\eps}{\delta}$. 
More precisely, since $C(x,r)$ is the same as a light cone with vertex at the point $(x,r)$, the bounds
in~\cite{S2} were obtained for families $\{(x_j,r_j)\}_{j=1}^N$ with $\delta$-separated $x_j$ and a unique~$r_j$ for every $x_j$. But the method was to estimate the {\em three-dimensional measure} of various sets of $(x,r)\in\R^3$ and then to project this bound onto the plane. The latter is always based on the fact that there
is a fixed ``slack'' in the vertical direction, which is precisely the variable corresponding to the radius.
For example, for the case of~\eqref{eq:int} it is clear that the amount of freedom in the set on the left hand side in the radial direction is $\eps$, which explains the factor of $\frac{\eps}{\delta}$. Hence the measure estimate in Lemma~2.5 of~\cite{S2}, which is 
\[ |\durch|\les\frac{ t^2}{\delta^2}\min(\sqrt{\frac{\epsilon}{\tau}},\frac{\epsilon}{\sqrt{\beta \tau}}) \]
for the {\em two-parameter situation}, needs to be multiplied by $\frac{\eps}{\delta}$, and we are done.
\eproof

\noindent The following lemma is our main technical lemma. Observe that \eqref{eq:hypo} holds
with $A=\delta^{-3}$, say. The desired Strichartz bound will then follow easily be iterating this lemma,
see Corollary~\ref{cor:main} below.

\begin{lemma}
\label{lem:iter}
Suppose with some constant $A\ge1$ (which may depend on $\delta$ but nothing else)
\begin{equation}
\label{eq:hypo}
|\{C^\delta\:|\: \mu^{\CCet}_\delta > A\lambda^{-1}|\C|^{\frac23}  \}| \le 
\lambda |C^\delta|
\end{equation}
for all $C\in\C$, $\delta\le \eps\le t$, $0<\lambda<1$. Then there exists $\A\subset\C$, $|\A|\ge\half|\C|$, 
so that
\begin{equation}
\label{eq:concl}
|\{C^\delta\:|\: \mu^{\Acet}_\delta > |\log\delta|^5\, C_0\, \sqrt{A} \,\lambda^{-1} |\A|^{\frac23} \}| \le \lambda|C^\delta|
\end{equation}
for all $C\in\A$, $\delta\le\eps\le t$, $0<\lambda<1$. Here $C_0$ is some absolute constant.
\end{lemma}
\noindent
\underline{Proof:} Let $N=|\C|$.  Suppose that
for at least half the circles in~$C\in \C$ one has (with $b=5$, say) 
\begin{equation}
\label{eq:large}
|\{C^\delta\:|\: \mu^{\CCet}_\delta > |\log\delta|^b C_0 \sqrt{A}\, \lambda^{-1} |\C|^{\frac23} \}| > \lambda|C^\delta|
\end{equation}
for some choice of $\eps, t$ and $\lambda$. Then pigeonhole to get $\B\subset\C$ and fixed $\delta\le\eps\les t\le 1$, $0<\lambda<1$ such that $|\B|\gtrsim |\log\delta|^{-3} |\C|$ and~\eqref{eq:large} 
holds for all $C\in\B$. We first claim that 
\begin{equation}
\label{eq:lambdabig} 
\lambda \ge \frac{\delta}{\sqrt{\eps t}}.
\end{equation}
The point here is that the lower bound is comparable to the area of intersection of $C_1^\delta$ and~$C_2^\delta$ with $\Delta=\eps$ and $d=t$, see~\eqref{eq:durch}.
 Indeed, for any fixed $x\in\R^2$ and $C\in\C$,  
\begin{equation}
\label{eq:triv} 
\mu^{\CCet}_\delta (x) \les  \min(\sqrt{t\eps}\;t\eps\,\delta^{-3},N) \les (\sqrt{t\eps}\: t \eps/\delta^3)^{\frac13}N^{\frac23}.
\end{equation}
The $\sqrt{t\eps}\;t\eps\,\delta^{-3}$ term derives from the fact that the centers of the circles have to lie
in a rectangles of dimensions~$\sqrt{t\eps}\,t$, whereas the radius has a freedom of~$\eps$. 
If $\lambda \le \frac{\delta}{\sqrt{\eps t}}$ and $x$ belongs to the set in~\eqref{eq:large}, then 
\[ \mu^{\CCet}_\delta (x) \ge |\log\delta|^b\,C_0\sqrt{A}\, \frac{\sqrt{\eps t}}{\delta} N^{\frac23},\]
which contradicts the apriori bound \eqref{eq:triv}. Hence~\eqref{eq:lambdabig} holds as claimed. This allows 
one to run a counting argument. More precisely, define
\bea
S &:=& \{(C,C_1,C_2)\:|\: C\in\B,\quad C_1,C_2\in \CCet, \quad\sgn(r-r_1)=\sgn(r-r_2), \nn \\
&& \qquad\qquad\qquad \beta\sim \Delta(C_1,C_2), \quad|x_1-x_2|\sim\tau\} \nn
\enea
where $C_1=C_1(x_1,r_1)$, $C_2=C_2(x_2,r_2)$, and $\beta,\tau$ are chosen by means of pigeonholing so that
\begin{equation}
\label{eq:stern} 
\card(S) \gtrsim |\log\delta|^{-2} |\B|\, \bigl(\mu_0\lambda\frac{\sqrt{\eps t}}{\delta}\bigr)^2
\gtrsim C_0^2\,|\log \delta|^{-5+2b} A\, N\,N^{\frac43}\, \frac{\eps t}{\delta^2},
\end{equation}
cf.~our assumption \eqref{eq:large}. Here 
\[ \mu_0 = C_0\,|\log\delta|^b\sqrt{A}\,\lambda^{-1}N^{\frac23}.\]
The reason for including the condition $\sgn(r-r_1)=\sgn(r-r_2)$ in the definition of the set $S$ is 
to ensure that for most $(C,C_1,C_2)\in S$ one has $C_1\cap C_2\not=\emptyset$. More precisely, it 
follows immediately from Lemma~2.6 in~\cite{S2} that
\bea
\card \big\{(C,C_1,C_2)\in S\:|\: C_1^\delta\cap C_2^\delta = \emptyset \big\} &\lesssim& N^2 
\frac{\sqrt{\eps t}\; t\eps}{\delta^3}
\nonumber \\
&\lesssim& N^2 (\sqrt{\eps t}/\delta)^2 N^{\frac13} \ll |\log\delta|^{-10}\,N (\mu_0\lambda \frac{\sqrt{\eps t}}{\delta})^2.
\nonumber
\enea
Hence
\[ S' := \{ (C,C_1,C_2)\in S\:|\:C_1^\delta\cap C_2^\delta \not=\emptyset, \;\Delta(C_1,C_2)\sim \beta, \; d(C_1,C_2)\sim \tau \}\]
has a lower bound on its cardinality that is comparable to the one in~\eqref{eq:stern}. For simplicity, we will henceforth write $S$ but mean $S'$.  
Suppose $\tau<\eps$. Then, since $\beta\les \tau+\eps\le 2\eps$, 
\[ \card\{C_2\in\C\:|\: \Delta(C_1,C_2)\sim\beta, |x_1-x_2|\sim \tau\} \les (\eps/\delta)^3.\]
In this case we bound the cardinality of $S$ from about by fixing $C,C_1$ and then choosing~$C_2$. 
This can be done in no more than 
$ N^2 \min\bigl(N, (\eps/\delta)^3\bigr) $
many ways. Thus,
\bea
|\log\delta|^{-5} N(\mu_0\lambda\frac{\sqrt{t\eps}}{\delta})^2 &\lesssim& \card(S) \les N^2 N^{\frac13} (\eps/\delta)^2 \qquad \Longrightarrow\nonumber \\
\mu_0 &\lesssim&  |\log \delta|^{\frac52} \lambda^{-1}N^{\frac23}\sqrt{\eps/t} \lesssim |\log\delta|^{\frac52} \lambda^{-1} N^{\frac23}, \label{eq:mu0}
\enea
which is a contradiction. Hence we may assume that $\tau\ge\eps$ so that $d(C_1,C_2)\sim |x_1-x_2|$.
Furthermore, suppose $N<\frac{(\eps t)^{\frac32}}{\delta^3}$. Then as in~\eqref{eq:mu0} one obtains
\bea
|\log\delta|^{-5} N(\mu_0\lambda\frac{\sqrt{t\eps}}{\delta})^2 &\lesssim& \card(S) \les N^2 N^{\frac13} 
\frac{\eps t}{\delta^2} \qquad \Longrightarrow\nonumber \\
\mu_0 &\lesssim& |\log\delta|^{\frac52} \lambda^{-1} \frac{\delta}{\sqrt{\eps t}} N^{\frac23} \frac{\sqrt{\eps t}}{\delta} 
= C_0\,|\log\delta|^{\frac52}\, \lambda^{-1}N^{\frac23}, \nn 
\enea
which is again a contradiction. 
So we can assume that $N>(\eps t)^{\frac32} \delta^{-3}$. 
We are now in a position to bound the cardinality of~$S$ from above. To do so, we
fix one of at most~$N$ choices of~$C_1$. 
It follows easily from our hypothesis~\eqref{eq:hypo} that
\[
\card\{C_2\in\C_{\beta \tau}^{C_1}\:|\: C_1\cap C_2\not=\emptyset\} \lesssim |\log\delta|\,A\,N^{\frac23}\frac{\sqrt{\beta\tau}}{\delta}.
\] 
Indeed, use \eqref{eq:durch}, \eqref{eq:hypo} and sum over $\delta\gtrsim\lambda=2^{-j}\les 1$.
This controls the number of choices of $C_2$. Finally, for a fixed pair $(C_1,C_2)$ we bound the number
of choices for~$C$ by means of  Lemma~\ref{lem:twocirc}. The details are as follows:

\smallskip\noindent
\underline{Case 1:} If $\beta<\eps$, then by Lemma~\ref{lem:twocirc}
\bea
\card (S) &\lesssim& N |\log\delta|A\,N^{\frac23} \frac{\sqrt{\beta\tau}}{\delta} \frac{\eps t^2}{\delta^3} \sqrt{\frac{\eps}{\tau}} \nonumber \\
&\le& N |\log\delta|A\, N^{\frac23} \frac{\eps^2 t^2}{\delta^4}  \label{eq:zwischen}\\
&\lesssim& N |\log\delta|A\, N^{\frac43} \frac{\eps t}{\delta^2}  \nonumber
\enea
where the final inequality uses our assumption $\eps t< \delta^2 N^{\frac23}$. 
This is a contradiction to~\eqref{eq:stern}.

\smallskip\noindent
\underline{Case 2:} If $\beta\ge \eps$, then again by Lemma~\ref{lem:twocirc}
\bea
\card(S) &\lesssim& N |\log\delta|A\, N^{\frac23} \frac{\sqrt{\beta\tau}}{\delta} \frac{\eps t^2}{\delta^3} \frac{\eps}{\sqrt{\beta\tau}} \nonumber \\
&\lesssim& N |\log\delta|A\, N^{\frac23} \frac{\eps^2 t^2}{\delta^4}, \nonumber
\enea
which is the same as \eqref{eq:zwischen}
and we are done. \eproof

\noindent Iteration of the previous lemma leads to the following corollary. It proves that
condition iii) of Lemma~\ref{lem:equiv} holds, and also thus the restricted weak type form
of Strichartz. 

\begin{cor}
\label{cor:main}
Given a $\delta$-separated three parameter family of circles $\C$, one has: For any $\eta>0$ and $\delta$ sufficiently small depending on~$\eta$, there exists $\A\subset\C$ such that $|\A|>\eta|\C|$ and
\begin{equation}
\label{eq:mult} 
|\{C^\delta\:|\: \mu_\delta^{\A} > \delta^{-4\eta}\lambda^{-1}|\A|^{\frac23}\}| \le \lambda|C^\delta|
\end{equation}
for all $C\in\A$, $0<\lambda<1$.
\end{cor}
\noindent
\underline{Proof:} Note that \eqref{eq:hypo} holds with $A=\delta^{-3}$. Applying Lemma~\ref{lem:iter} 
repeatedly, say~$K$ times, produces a subset $\A\subset\C$ of cardinality at least $2^{-K}|\C|$ so that 
\be
\nn
|\{C^\delta\:|\: \mu^{\Acet}_\delta > |\log\delta|^{10}\, C_0^2\, \delta^{-3.2^{-K}} \,\lambda^{-1} |\A|^{\frac23} \}| \le \lambda|C^\delta|
\end{equation}
for all $C\in\A$, $\delta\le\eps\le t$, $0<\lambda<1$. Hence
\bea
&& |\{C^\delta\:|\: \mu^{\A}_\delta > |\log\delta|^{12}\, C_0^2\, \delta^{-3.2^{-K}} \,\lambda^{-1} |\A|^{\frac23} \} | \nn \\
&&  \les \sum_{\eps,t} |\{C^\delta\:|\: \mu^{\Acet}_\delta > |\log\delta|^{10}\, C_0^2\, \delta^{-3.2^{-K}} \,\lambda^{-1} |\A|^{\frac23} \}| \nn \\
&& \les |\log\delta|^2\, \lambda |C^\delta|, \nn
\enea
for all $C\in\A$, $0<\lambda<1$. The sum here is over dyadic $\eps,t$. Absorbing the $|\log\delta|^2$-factor into $\lambda$ and taking $K$ sufficiently large finishes the proof.
\eproof


\section{A simplified proof of Wolff's $L^3_r(L^\infty_x)$ bound}
\label{sec:L3}

\noindent The bound considered in the previous section is only one out of many estimates dealing
with circular averages.
More precisely, one can ask about exponents $p,q,s$ so that
\bea
 \|\sigma_{rS^1}\ast f\|_{L_r^p(L_x^q)} &\les& \|f\|_{L^s(\R^2)}  \label{eq:xfirst} \\
 \|\sigma_{rS^1}\ast f\|_{L_x^q(L_r^p)} &\les& \|f\|_{L^s(\R^2)},  \label{eq:rfirst}
\enea
where $\sigma_{rS^1}$ is the normalized measure on the circle $rS^1$ of radius~$r$.
The correct range of exponents for these bounds can be found by means of the usual examples, namely,
the focusing, Knapp, and scaling examples. These refer to setting $f=\chi_{C^\delta(0,1)}$, $f=\chi_R$ where
$R$ is a $\delta\times\sqrt{\delta}$-rectangle, and $f=\chi_{B(0,\delta)}$, respectively.
We will completely ignore the second class of estimates, i.e., \eqref{eq:rfirst}, in this paper.
Let it suffice to say that the endpoint $p=\infty$, $s=q>2$ is Bourgain's circular maximal theorem, 
see~\cite{B} and~\cite{MSS}. As far as the first class~\eqref{eq:xfirst} is concerned, it is easy to check that one endpoint is the
estimate~\eqref{eq:3/2}, i.e., $p=\infty, q=3, s=\frac32$. The other is the case $q=\infty$, $p=s=3$. 
In fact, in \cite{W} Wolff proved that for every $\eps>0$ 
\begin{equation}
\label{eq:L3}
\|\sup_{x\in \R^2} (\sigma_{rS^1} \ast f)(x) \|_{L^3_r([1,2])} \le C_\eps\, \|f\|_{W^{3,\eps}(\R^2)}
\end{equation}
by means of a combinatorial device originating in~\cite{cells} called the method of cell decomposition. 
It is not our intention to review this method, as the paper~\cite{cells} is highly readable, and 
because Wolff explains his adaptation of it in~\cite{W} and to lesser extent also in~\cite{tomrev} and~\cite{W2}. 
Note that the $\eps$ is necessary in~\eqref{eq:L3} since there exist sets of measure zero that contain
a circle of every radius, see~\cite{tomrev}. \\
The paper~\cite{W} is rather complicated, but~\cite{W2} allows for some significant simplifications. 
It is pointed out in~\cite{W2} that the main estimate from Section~1 of that paper allows for a 
simplified proof of~\eqref{eq:L3} be means of ``fairly standard arguments''. While this is true on a heuristic
level, it is perhaps less true on the level of a rigorous argument. We therefore hope that the
proof of this fact presented here is of some value. 

\noindent We start by recalling some terminology from~\cite{W2}: 
\begin{itemize}
\item
Let $\W$ and $\B$ be families
of circles which are each $\delta$-separated. We refer to the pair $\W,\B$ as $t$-bipartite
provided $t\le d(w,b)\le  100t$ if $w\in\W$ and
$b\in\B$, $d(w_1,w_2)\le t$ for $w_1,w_2\in\W$,  and $d(b_1,b_2)\le t$ for $b_1,b_2\in\B$. 
\item
A $(\delta,t)$-rectangle is a $\delta$-neighborhood of an arc of length~$\sqrt{\frac{\delta}{t}}$ on
some circle. It follows from~\eqref{eq:durch} that two annuli $C^\delta_1,C^\delta_2$ with 
$\Delta(C_1,C_2)\le \delta$ and $t\le d(C_1,C_2)\le 2t$ intersect in a set that can be
covered by a finite number (some absolute constant) of $(\delta,t)$-rectangles.
\item We say that two $(\delta,t)$-rectangles are comparable, if there is an $(a_0\delta,t)$-rectangle that
contains them both where $a_0$ is some absolute constant. A circle $C$ is said to be tangent
to a $(\delta,t)$-rectangle~$R$ if the $a_1\delta$-neighborhood of $C$ contains~$R$, where $a_1$ is some
fixed constant. A $(\delta,t)$-rectangle
$R$ is said to be of type $(\ge\mu,\ge\nu)$ relative to a $t$-bipartite pair $\W,\B$ as above provided
there are at least $\mu$ circles from~$\W$ and at least~$\nu$ circles from $\B$ that are tangent to~$R$.
\end{itemize}

\noindent We refer the reader to Section~1 of~\cite{W2} for more details. We will need the following
estimate which is Lemma~1.4 in that paper. 

\smallskip\noindent
\underline{Bound on high multiplicity rectangles: } Let $\W,\B$ be a $t$-bipartite pair. If $\eps>0$ then there is a constant $C_{\eps}$ such that the cardinality
of any set of pairwise incomparable $(\delta, t)$-rectangles of type $(\geq\mu,\geq\nu)$ relative to $\W,\B$ 
is bounded by 
\begin{equation}
\label{eq:wow} 
C_{\eps}(mn)^{\eps}\left( (\frac{mn}{\mu\nu})^{\frac{3}{4}}+\frac{m}{\mu}+\frac{n}{\nu}\right)
\end{equation}
where $m=|\W|$ and $n=|\B|$.

\smallskip 
It is evident that~\eqref{eq:wow} allows one to control the number of pairs $(w,b)\in \W\times\B$ which
are $\delta$-tangent (which means that $\Delta(w,b)\les\delta$). Indeed, by the second item above, counting
pairs of $\delta$-tangent circles (i.e., counting incidences) 
is the same as counting incomparable $(\delta,t)$-rectangles which are obtained as intersections
of at least one $\delta$-annulus from~$\W$ with another from~$\W$. Of course, one has to keep track
of multiplicity here. For example, if no two circles in $\W$ or $\B$ are $\delta$-tangent, then~\eqref{eq:wow}
with $\mu=\nu=1$ gives
\[ \#\{(w,b)\in \W\times\B\:|\: \Delta(w,b)\le \delta,\; t\le d(w,b)\le 2t\} \les (|m||n|)^{\frac34+\eps}.\]
This is the analogue of a bound which is implicit in~\cite{cells} (see also~\cite{tomrev}) 
that says that for any collection of $N$ circles in
the plane so that no three are tangent at the same point the total number pairs of {\em exactly} tangent 
circles is at most $N^{\frac32+\eps}$. If, on the other hand, $\mu=m$ and $\nu=n$ which is the case where all are tangent to a single rectangle, then the bound 
in~\eqref{eq:wow} is~$\delta^{-\eps}$. \\
We now start with the proof that~\eqref{eq:wow} implies~\eqref{eq:L3}. 
Firstly, the same technique that was used in Lemma~\ref{lem:equiv} shows that~\eqref{eq:L3} 
is equivalent with the statement of  the following lemma, 
which is the ``main lemma'' in Wolff~\cite{W2}, see page~998.

\begin{lemma}
\label{lem:reduc}
Given $\eta>0$ the following holds for sufficiently small $\delta$:
Suppose $\C$ is a family of circles with $\delta$-separated radii. Then there exists
$\A\subset\C$, $|\A|>\half |\C|$ such that
\begin{equation}
\label{eq:L3mult}
|\{C^\delta\:|\: \mu^{\C}_\delta > \delta^{-\eta} \lambda^{-2}\}| \le \lambda|C^\delta|
\end{equation}
for all $C\in\A$ and $0<\lambda<1$.
\end{lemma}
\noindent
\underline{Proof using \eqref{eq:wow}:} Heuristically speaking, this is very simple. Indeed,
suppose that a typical annulus $C^\delta$ from~$\C$ contains about $A$ many $\delta\times\sqrt{\delta}$-rectangles  each of which has about $\lambda^{-2}$ many circles tangent to it (assuming $t\sim1$ here). Then~\eqref{eq:wow} implies the following bound on the number of pairs of tangent circles, with $N=|\C|$ and $\mu=\lambda^{-2}$: 
\[ N\,A\lambda^{-2} \les (\frac{N}{\mu})^{\frac32} \mu^2 = N^{\frac32} \lambda^{-1}.\]
Hence $A\les N^{\frac12} \lambda$ and since  $N\les \delta^{-1}$ one obtains
\[ |\{C^\delta\:|\: \mu^{\C}_\delta > \lambda^{-2}\}| \les A\delta^{\frac32} \les \sqrt{N}\lambda\,\delta^{\frac32}\les \lambda\,\delta.\]

\noindent To make this argument rigorous, we shall use induction in~$\delta$. I.e., let $\delta>0$ be small and
assume that the statement holds for~$2^{j}\delta$, $j\ge1$ (the case of $\delta\sim1$ being trivial).  
Now suppose that at least half the circles $C\in\C$ satisfy
\begin{equation}
\label{eq:first}
 |\{C^\delta\:|\: \mu^{\C}_\delta > \delta^{-\eta}\lambda^{-2}\}| > \lambda|C^\delta|
\end{equation}
for some $0<\lambda<1$ (depending on~$C$). 
Then there exist fixed choices of $\delta\le\eps\les t \le1$ and $0<\lambda_0<1$
such that
\begin{equation}
\label{eq:eta/2} 
|\{x\in C^\delta\:|\: \mu^{\CCet}_\delta(x) \sim \mu_0,\quad \mu_\delta^{\C}(x)\les \mu_0|\log\delta|^2 \}| > \lambda_0|C^\delta|
\end{equation}
for all $C\in\A$ where 
\begin{equation}
\label{eq:Acard} 
|\A|>|\log\delta|^{-3}|\C| \text{\ \ and\ \ }\mu_0\gtrsim \delta^{-\frac{\eta}{2}}\lambda_0^{-2}.
\end{equation}
Assume first that $\delta\le\eps<2\delta$. In what follows let
\[ \B^{C}_\delta := \bigcup_{s}\B_{\delta s}^{C} \qquad \W^{C}_\delta := \bigcup_{s} \W_{\delta s}^{C}.\] 
We now show that there exist $\W,\B\subset\C$ so that
\begin{enumerate}
\item $\W,\B$ is a $\frac{t}{4}$-bipartite pair
\item for all $C\in\W$ one has
\begin{equation}
\label{eq:vielf}
|\{x\in C^\delta\:|\: \mu_\delta^{\B^C_\delta}(x)\sim \mu_2,\; \mu_\delta^{\W^C_\delta}(x)\les \mu_1\}| > \lambda_0|C^\delta|
\end{equation}
where $\mu_2\gtrsim \delta^{-\frac{\eta}{2}}\lambda_0^{-2}$ and $1\le \mu_1\les |\log\delta|^2\mu_2$. 
\item $|\W|\gtrsim |\log\delta|^{-4}|\B|>>1$.
\end{enumerate} 
Strictly speaking, in ii) one needs to write $\lambda_0|\log\delta|^{-1}$ on the right--hand side, but we can
absorb the $|\log\delta|$-factor into~$\lambda_0$. Let $\A=\{C(x_j,r_j)\}_{j=1}^M$. 
Then set
\bea 
 \W &=& \{C(x_j,r_j)\in\A\:|\: (x_j,r_j)\in Q\} \nonumber \\
 \B &=& \{C(x_j,r_j)\in\C\:|\: \frac{11t}{20}<\dist((x_j,r_j),Q) < \frac{19t}{20}\} \label{eq:WBdef}
\enea
where $Q$ is a ball of size $\frac{t}{20}$ in $\R^3$ for which iii) holds. To see that such a ball exists,
consider a covering on $\R^3$ by balls $Q$ of this size which have overlap bounded by some absolute constant.
Simultaneously, consider a covering by balls~$Q^*$ with the same centers as the $Q$'s but ten times their size.
In view of~\eqref{eq:Acard} there has to exist one such ball~$Q$ so that 
\begin{equation}
\label{eq:locdens} 
|\A\cap Q| \gtrsim |\log\delta|^{-3} |\C\cap Q^*|
\end{equation}
which implies iii) with $|\log\delta|^{-3}$. 
By construction, if $C\in\W$ and $\Cb\in\B^C_\delta$, then $\frac{t}{2}\le d(C,\Cb)\le t$. 
Thus~\eqref{eq:eta/2} implies that
\begin{equation}
\label{eq:againmult} 
|\{C^\delta\:|\: \mu_\delta^{\B^C_\delta} \gtrsim \mu_0, \; \mu_\delta^{\W}\les \mu_0|\log\delta|^2\}| > \lambda_0|C^\delta|,
\end{equation}
see \eqref{eq:CCetdef}. 
Pigeonholing again one obtains~\eqref{eq:vielf} at the cost of replacing $\lambda_0$ with $\lambda_0|\log\delta|^{-1}$, see the comment above, as well as with a loss of another factor of $|\log\delta|^{-1}$ in~iii).
 Finally, satisfying~i) requires one more application of the pigeonhole principle,
but only with finitely many cases. Indeed, cover $\B$ with a finite number (some absolute constant) of balls
of size~$\frac{t}{4}$ and replace $\B$ with the intersection of itself with one of these balls for
which~ii) remains true.

\noindent We now distinguish two cases, namely $\lambda_0\gg\sqrt{\frac{\delta}{t}}$ and $\lambda_0\les\sqrt{\frac{\delta}{t}}$. 

\smallskip\noindent
\underline{Case 1:} $\lambda_0 \gg \sqrt{\frac{\delta}{t}}$

\noindent In this case we count pairs of $\delta$-tangent circles $(C_1,C_2)\in \W\times\B$.
By~\eqref{eq:vielf} every $C^\delta$ with $C\in\W$ contains at least $A$ many $(\delta,t)$-rectangles, 
\begin{equation}
\label{eq:Aatleast}
A \gtrsim \lambda_0 \sqrt{\frac{t}{\delta}},
\end{equation}
each of which is $\delta$-tangent to about $\mu_2$ many circles
from~$\B$ and to at most about~$\mu_1$ many circles from~$\W$ (but clearly, tangent to at least one circle from~$\W$). Thus, from~\eqref{eq:wow}, 
\bea
|\W|\,A\,\mu_2 &\les& \delta^{-\frac{\eta}{100}}\sum_{2^j\les\mu_1} \left(\frac{|\W||\B|}{2^j\,\mu_2}\right)^{\frac34} 2^j\,\mu_2 \nonumber \\
  &\les& \delta^{-\frac{\eta}{90}} |\W|^{\frac32}(\mu_1\mu_2)^{\frac14} \les \delta^{-\frac{\eta}{80}}
 |\W|^{\frac32} \mu_2^{\half}.  \nonumber
\enea
Hence, since $|\W|<t\delta^{-1}$, 
\[ A \les \delta^{-\frac{\eta}{80}}|\W|^{\frac12}\mu_2^{-\half} \les \delta^{-\half+\frac{\eta}{5}} \lambda_0\sqrt{t}, \]
which contradicts \eqref{eq:Aatleast}. 

\smallskip\noindent
\underline{Case 2:} $\lambda_0\les\sqrt{\frac{\delta}{t}}$

\noindent Fix any $C\in\W$ and let $x$ be a point from the left-hand side of~\eqref{eq:vielf}.
Then 
\[ 
\frac{t}{\delta}\gtrsim |\B|\ge \mu_\delta^{\B_\delta^{C}}(x) \gtrsim \mu_2 
\gtrsim \delta^{-\frac{\eta}{2}} \frac{t}{\delta},
\]
which is impossible. Hence we are done with the case $\delta\le \eps\le 2\delta$. 

\noindent To treat the case $\eps\ge2\delta$ we will need to ``thin'' the sets $\W,\B$ at scale~$\eps$
in order to apply~\eqref{eq:wow} to $\eps$-tangencies. 
More precisely, define $\W$ and~$\B$ as in~\eqref{eq:WBdef}, \eqref{eq:locdens}. 
In particular, \eqref{eq:againmult} holds. It will be convenient to pass to subsets of $\W,\B$ 
that are homogeneous at scale~$\eps$ in the radial variable.
Partition $\R^3$ into disjoint slabs of size $\eps$, i.e., $\R^3=\bigcup_{\ell\in\Z} S_\ell$
where 
\[ S_\ell = \{(x,t)\in\R^3\:|\: \ell\eps < t \le (\ell+1)\eps \}.\]
For every $S_\ell$ one has $\card(S_\ell\cap\W) \sim \rho_\ell$ for some 
\begin{equation}
\label{eq:1para} 
1\le\rho_\ell\les \frac{\eps}{\delta}
\end{equation}
provided the intersection is not empty (since the circles in $\W$ have $\delta$-separated radii). 
Thus there exists $\tau_{\W}$ such that 
\[ \Bigl|\;{\bigcup}'_{\ell:\rho_\ell\tau_{\W} \sim 1} S_\ell\cap\W \;\Bigr| \gtrsim |\log\delta|^{-1}|\W|\]
where $\bigcup'$ means that $\ell$ is required to be either even or odd, depending on which
choice leads to the larger set. This ensures that the separation between points in different
slabs is bigger than~$\eps$. 
Denote the set on the left--hand side by $\W_{\rm hom}$. Similarly, there exists $\tau_{\B}$ 
so that~\eqref{eq:vielf} 
remains correct (up to logarithmic factors) if $\B$ is replaced with~$\B_{\rm hom}$, the latter being
\[ {\bigcup}'_{\ell:\rho_\ell\tau_{\B} \sim 1} S_\ell\cap\B. \]
Notice that~iii) above only changes by another $|\log\delta|$-factor. For simplicity, we will
ignore logarithmic factors altogether from now on. Moreover, in view of the preceding, conditions
i), ii), iii) above remain valid with a suitable choice of $\mu_1,\mu_2$ if we 
replace $\W$ with~$\W_{\rm hom}$ and $\B$ with~$\B_{\rm hom}$, and we will drop the ``hom'' from now on. 
Define $\wt{\W}\subset\W$ by randomly selecting one point from each nonempty
$S_\ell\cap\W$, and similarly~$\wt{\B}\subset\B$. By the homogeneity property,
\begin{equation}
\label{eq:cardthin}
|\wt{\W}|\sim \tau_{\W}|\W| \text{\ \ and\ \ } |\wt{\B}| \sim \tau_{\B}|\B|.
\end{equation}
This holds for {\em every} choice of points in $\wt{\W}$ and $\wt{\B}$ --- the reason for choosing the
points randomly rather than deterministically will become clear only later on. 
As before, we will count pairs $(C_1,C_2)\in \wt{\W}\times \wt{B}$ that are $\eps$-tangent 
by means of~\eqref{eq:wow}. Such pairs meet in $(\eps,t)$-rectangles of which a typical
one is tangent to about $\wt{\mu_2}$ circles from~$\B$, and no more than~$\wt{\mu_1}$ circles from~$\W$,
respectively. We will need to determine~$\wt{\mu_1}$ and~$\wt{\mu_2}$; more precisely, we will 
bound $\wt{\mu_1}$ from above, and $\wt{\mu_2}$ from below. 
By the induction hypothesis applied to $\wt{\W}$ (a set of circles with $\eps$-separated
radii),  at least half the circle $C\in\wt{\W}$ satisfy (with a fixed choice of small $\eta>0$) 
\begin{equation}
\label{eq:goodW} 
|\{C^{\eps}\:|\: \mu^{\wt{\W}}_\eps > \eps^{-\frac{\eta}{100}} \lambda^{-2}\}| \le \lambda|C^\eps|
\end{equation}
for all $0<\lambda<1$. We now replace $\wt{W}$ with this subcollection of circles; for convenience, we
again denote it by~$\wt{\W}$, which then satisfies~\eqref{eq:goodW} afortiori. 
In what follows, $R$ denotes an $(\eps,t)$-rectangle which is
$\eps$-tangent to a pair $(C_1,C_2)\in\wt{\W}\times\wt{\B}$. Note that on such a rectangle the functions
\[ \mu^{\wt{\W}^{C_1}_\eps}_\eps \text{\ \ \ and\ \ \ } \mu^{\wt{\B}^{C_2}_\eps}_\eps \]
are basically constant. 

\smallskip\noindent
\underline{Case 1:} $\lambda_0>>\sqrt{\frac{\eps}{t}}$. 

\noindent Define
\[ {\mathcal S} := \{(C_1,C_2)\in \wt{\W}\times\wt{\B} \:|\: C_1,C_2 \;\text{\ are $\eps$-tangent to~$R$ and\ } \mu^{\wt{\W}^{C_1}_\eps}_\eps \le \wt{\mu_1},\; \mu^{\wt{\B}^{C_2}_\eps}_\eps\sim \wt{\mu_2}\;\text{\ on $R$}\}\]
where 
\begin{equation}
\label{eq:mu1mu2}
\wt{\mu_1} :=\delta^{-\frac{\eta}{50}}\lambda_0^{-2}\text{\ \ and\ \ }\wt{\mu_2}\gtrsim\frac{\mu_2}{\nu}\frac{\eps}{\delta}\tau_{\B}
\end{equation}
for some positive integer $\nu$. To understand these values, fix any $C\in\wt{\W}$. From~ii) above and another application of the pigeonhole principle one concludes that there exist 
\begin{equation}
\label{eq:Anu} 
A\ge \nu \lambda_0 \sqrt{\frac{t}{\eps}}
\end{equation}
many $(\eps,t)$-rectangles $R$ with the property that each one of them is hit by about 
\begin{equation}
\label{eq:duenn} 
\frac{\mu_2}{\nu} \frac{\delta \sqrt{\frac{\eps}{t}}}{\delta^2/\sqrt{\eps t}} = \frac{\mu_2}{\nu} \frac{\eps}{\delta} 
\end{equation}
many annuli in $\B^{C}_{\delta}$. Denote the collection of these $(\eps,t)$-rectangles by~${\mathcal R}(C)$.
Recall that the set $\wt{\B}$ is defined by random selection of points from~$\B$, see~\eqref{eq:cardthin}. 
In what follows, ``probability'' refers to this random selection.

\noindent We now claim the following: With high probability, at least half the circles $C\in\wt{\W}$ have the property that all the rectangles in~${\mathcal R}(C)$ are hit by 
\begin{equation}
\label{eq:uff} 
  \gtrsim \frac{\mu_2}{\nu}\frac{\eps}{\delta}\tau_{\B} 
\end{equation}
annuli in $\wt{\B}^{C}_{\eps}$ (which explains the choice of $\wt{\mu_2}$ in~\eqref{eq:mu1mu2}).
This follows from an elementary large deviation estimate for Bernoulli variables, 
see Lemma~\ref{lem:bern} below. Firstly, observe
that~\eqref{eq:Anu} and~$A\lesssim \sqrt{\frac{t}{\eps}}$ imply
\[ \nu \lesssim \lambda_0^{-1}.\]
Since $\tau_{\B}\gtrsim \delta/\eps$,  the right-hand side of~\eqref{eq:uff} is  $\gtrsim \mu_2/\nu \gtrsim \delta^{-\frac{\eta}{2}}$. If $R\in{\mathcal R}(C)$ is fixed, then
we say that $R$~{\em is good}, provided the number of circles from~$\wt{\B}^C_\eps$ that intersect~$R$ is
at least~$\delta^{\frac{\eta}{100}}$ times the expected number of circles hitting it. But in view of~\eqref{eq:duenn} the latter is at least
\[ \frac{\mu_2}{\nu}\frac{\eps}{\delta}\tau_{\B} \gtrsim \delta^{-\frac{\eta}{2}}.\]
Thus a good rectangle is hit by at least
\[ \delta^{\frac{\eta}{100}}\delta^{-\frac{\eta}{2}} \ge \delta^{-\frac{\eta}{3}}\]
many circles from~$\wt{\B}^C_\eps$. It follows from Lemma~\ref{lem:bern} that for every $R\in{\mathcal R}(C)$
\[  \Prob[R \text{\ \ is bad}] \les \exp(-c\,\delta^{-\frac{\eta}{3}}).\]
Thus also, for every $C\in\wt{W}$, 
\[ \Prob[\text{there exists a bad $R\in{\mathcal R}(C)$}] \les \exp(-c\,\delta^{-\frac{\eta}{3}}),\]
and so the claim under~\eqref{eq:uff} holds.
For any $C\in\wt{\W}$ it follows from~\eqref{eq:goodW} above that at most half  of the rectangles $R\in{\mathcal R}(C)$ satisfy $\max_{x\in R}\mu_\eps^{\wt{\W}^C_\eps}(x)>\eps^{-\frac{\eta}{100}}\lambda_0^{-2}$. In other words, at least $\half A$ many
rectangles are $\eps$-tangent to at most~$\eps^{-\frac{\eta}{100}}\lambda_0^{-2}$ circles from~$\wt{\W}$. 
Therefore,
\bea
\card({\mathcal S}) &\gtrsim& |\wt{\W}|\,A\,\wt{\mu_2} \sim \tau_{\W}|\W|\,A\,\wt{\mu_2}\label{eq:Slower}  \\
\card({\mathcal S}) &\les& \delta^{-\frac{\eta}{100}} \left(\frac{|\wt{\W}||\wt{\B}|}{\wt{\mu_1}\wt{\mu_2}}\right)^{\frac34}\, \wt{\mu_1}\wt{\mu_2} \nonumber \\
&\les& \delta^{-\frac{\eta}{100}} \left(\frac{\tau_{\W}|\W||\B|}{\frac{\mu_2}{\nu}\frac{\eps}{\delta}}\right)^{\frac34}\, (\delta^{-\frac{\eta}{50}}\lambda_0^{-2})^{\frac14}\wt{\mu_2}. \label{eq:Supper}
\enea
Recall that $\mu_2 \gtrsim \delta^{-\frac{\eta}{2}}\lambda_0^{-2}$ , see ii) above. 
Using that $\frac{t}{\delta}\gtrsim|{\W}|\gtrsim |{\B}|$ (see iii) above) as well as $\frac{\delta}{\eps}\les \tau_{\W} \les 1$
(see \eqref{eq:1para}), comparison of~\eqref{eq:Supper} with~\eqref{eq:Slower} thus leads to
\bea
 A &\les& \delta^{-\frac{\eta}{50}} \tau_{\W}^{-\frac14}|\W|^{\frac12} \lambda_0^{-\half}\Bigl(\frac{\mu_2}{\nu}\frac{\eps}{\delta}\Bigr)^{-\frac34} \nonumber \\
&\les& \delta^{\frac{\eta}{4}} \Bigl(\frac{\eps}{\delta}\Bigr)^{\frac14} \Bigl(\frac{t}{\delta}\Bigr)^{\half}
\lambda_0^{-\half}\lambda_0^{\frac32} \nu^{\frac34} \Bigl(\frac{\delta}{\eps}\Bigr)^{\frac34} \nonumber \\
&\les& \delta^{\frac{\eta}{4}} \Bigl(\frac{t}{\eps}\Bigr)^{\half} \lambda_0\,\nu^{\frac34},
\enea
which contradicts \eqref{eq:Anu}.

\smallskip\noindent
\underline{Case 2:} $\lambda_0 \les \sqrt{\frac{\eps}{t}}$. 

\noindent In this case there are at least
\[ A\ge \nu \ge 1\]
many $(\eps,t)$-rectangles with the property that each one of them is hit by
\bea &\gtrsim& \frac{\mu_2}{\nu}\frac{\lambda_0\delta}{\delta^2/\sqrt{\eps t}} = \frac{\mu_2}{\nu}\,\lambda_0\,
  \frac{\sqrt{\eps t}}{\delta} \nonumber \\
 &\gtrsim& \delta^{-\frac{\eta}{2}}\lambda_0^{-1}\,\frac{\sqrt{\eps t}}{\nu\delta}\gtrsim \delta^{-\frac{\eta}{2}} \frac{t}{\nu \delta}. \nonumber
\enea
many annuli in $\B_\delta^C$. But then,
\[ \frac{t}{\delta} \gtrsim |\B^C_\delta| \gtrsim \nu \delta^{-\frac{\eta}{2}} \frac{t}{\nu \delta},\]
a contradiction.
\eproof

\noindent The following lemma is a standard large deviation estimate for Bernoulli
variables.

\begin{lemma}
\label{lem:bern}
Let $X_1,X_2,\ldots,X_N$ be independent Bernoulli with 
\[ \Prob[X_1=1]=p \text{\ \ and\ \ } \Prob[X_1=0]=1-p\] 
for some $0\le p\le \half$. Then there exist absolute constants $c,C$ so that
\begin{equation}
\label{eq:ldt}
\Prob\bigl[\sum_{j=1}^N X_j < \alpha Np\bigr] \le C\exp(-c\, Np)
\end{equation}
for any $0<\alpha\le\half$. 
\end{lemma}
\underline{Proof: } Using Stirling, the probability in \eqref{eq:ldt} is estimated by
\bea
&& \sum_{\ell<\alpha Np} \frac{N!}{\ell!\,(N-\ell)!} p^\ell (1-p)^{N-\ell} \nn \\
&\les& 
\sum_{\ell<\alpha Np} \bigl(\frac{\ell}{N}\bigr)^{-\ell} \bigl(1-\frac{\ell}{N}\bigr)^{-(N-\ell)}
\exp\bigl(N\bigr[\frac{\ell}{N}\log p + \frac{N-\ell}{N}\log(1-p)\bigr]\bigr) \nn \\
& \les& \sum_{\ell<\alpha Np} \exp\bigl(N[H(\ell/N)-g_p(\ell/N)]\bigr) \label{eq:Hgp}
\enea
where $H(x)=-x\log x - (1-x)\log(1-x)$ and $g_p(x)=-x\log p -(1-x)\log(1-p)$.
Let $\phi_p(x)=H(x)-g_p(x)$. Then $\phi_p$ is increasing  on $[0,p]$, $\phi_p(p)=0$, $\phi_p'(p)=0$, and
$\phi''(p)=-\frac{1}{p(1-p)}$. Finally, one checks that $\phi'''(x)=\frac{1-2x}{x^2(1-x)^2}>0$ if $0<x<\half$.
It follows that $\phi''(x)\le \phi''(p)$ for $0\le x\le p$ so that
\[ \phi_p(x) \le \half (x-p)^2 \phi''_p(p) = - \frac{(x-p)^2}{2p(1-p)}\]
for those $x$. In particular, if $0\le x\le \alpha p$, then
\[ \phi_p(x) \le  - \frac{(1-\alpha)^2 p}{2(1-p)}.\]
Hence, the expression in \eqref{eq:Hgp} is no larger than
\[ N\exp\Bigl(-pN \frac{(1-\alpha)^2}{2(1-p)} \Bigr) \les \exp(-c\,pN),\]
since $\alpha\le\half$ and $p\le\half$. The lemma follows.
\eproof

\section{A lower bound for the tangency problem}

\noindent Wolff remarks that it is an interesting problem to decide whether the exponent~$\frac34$
in~\eqref{eq:wow} can be replaced by a smaller one. He also suggests that~$\frac23$ seems the most
reasonable conjecture for a sharp bound. This exponent is the same as the conjectured optimal
exponent for the Erd\"os unit distance problem in~$\R^3$, see~\cite{cells}, \cite{tomrev}. 
In this section we present examples that show why $\frac43$ would be optimal for the tangency problem (the distinction between
$\frac23$ and $\frac43$ comes from bipartite vs.~nonbipartite). 
These examples are simple and most likely standard, see~\cite{tomrev}, 
but we present them nevertheless. We start with
the exact tangency problem, i.e., given a collection of circles in the plane, decide how many pairs
of tangent circles there can be at most. Needless to say, we are interested in tangencies at those
points where only a small number of circles are tangent (low multiplicity).

Consider lattice points $(x,y,r)\in [1,N]^2\times [N,2N]$. The distinct circles
$C(x,y,r)$ and $C(x',y',r')$ are tangent iff
\[ |(x,y)-(x',y')| = |r-r'|. \]
So $(a,b,c):=(x-x',y-y',r-r')$ is an integer vector such that
\begin{equation}
\label{eq:tang} 
|a|,|b|,|c|\le N \text{\ \ and\ \ } a^2+b^2=c^2 \not=0.
\end{equation}
The following lemma is a well-known representation of Pythagorean triples.

\begin{lemma} 
\label{lem:pythago}
All integer solutions $(a,b,c)$ to  
\begin{equation}
\label{eq:pythtrip}
a^2+b^2=c^2,\quad b,c>0,\quad \gcd(|a|,b,c)=1
\end{equation}
are given by
\bea
\label{eq:sol}
&& a=\half(\alpha^2-\beta^2),\quad b=\alpha\beta,\quad c=\half(\alpha^2+\beta^2) \\
&& a=2^j\,\alpha^2-\beta^2,\quad b = 2^{\frac{j+2}{2}}\,\alpha\beta,\quad c=2^j\,\alpha^2+\beta^2, 
 \nn\\
&& a=\alpha^2-2^j\,\beta^2,\quad b = 2^{\frac{j+2}{2}}\,\alpha\beta,\quad c=\alpha^2+2^j\,\beta^2, 
 \nn
\enea
where $\alpha,\beta >0$  are odd integers, $\gcd(\alpha,\beta)=1$, and $j\ge2$ is even.
Moreover, this representation is unique.
\end{lemma}
\noindent
\underline{Proof:} Let $(a,b,c)$ be as stated. Then $b^2=(c-a)(c+a)$, and $c-a\not=0$, $c+a\not=0$.
Suppose $p>2$ is prime and $p|c-a$. Then $p^2|b^2=(c-a)(c+a)$. If $p|c+a$, then also $p|a$ and $p|c$
contradicting $\gcd(|a|,b,c)=1$. Hence $p^2|c-a$. In other words, 
$c-a=2^k\,\beta^2$ for some odd integer~$\beta>0$ and $k\ge0$.
Similarly, $c+a=2^\ell\,\alpha^2$ for some odd integer~$\alpha>0$ and $\ell\ge0$.  It follows that 
\[ 2a = 2^\ell\alpha^2-2^k\beta^2,\qquad 2c=2^\ell\alpha^2+2^k\beta^2\]
and thus
\[ b^2 = 2^{k+\ell}\,\alpha^2\beta^2.\]
Clearly, one needs $2|k+\ell$, and $\min(k,\ell)\le1$ (otherwise $2|a,b,c$). If $\min(k,\ell)=0$, then
in fact $k=\ell=0$ (otherwise $a,b$ are not integers). Thus,
\begin{equation}
\label{eq:rep1} 
a=\half(\alpha^2-\beta^2),\quad b=\alpha\beta,\quad c=\half(\alpha^2+\beta^2).
\end{equation}
In order to assure that $a,b,c$ are relatively prime, one needs $\gcd(\alpha,\beta)=1$. Conversely,
for any such $\alpha,\beta$ it follows that the representation~\eqref{eq:rep1} indeed gives a solution
of~\eqref{eq:pythtrip}. 
If $\min(k,\ell)=1$, then 
\begin{equation}
\label{eq:rep2} 
a=2^{\ell-1}\alpha^2-2^{k-1}\beta^2,\quad c=2^{\ell-1}\alpha^2+2^{k-1}\beta^2,\quad b=2^{\frac{k+\ell}{2}}\alpha\beta.
\end{equation}
If $k=\ell=1$, then $a,b,c$ would be even which is impossible. On the other hand, if $k\not=\ell$, then
$a,b$ are odd whereas $b$ is even. The necessary condition $\gcd(\alpha,\beta)=1$ remains valid in this case
as well. Conversely, under this condition and the restrictions $\min(k,\ell)=1$, $k\not=\ell$, $k+\ell$ even,
one checks that $\gcd(|a|,b,c)=1$ and that~\eqref{eq:rep2} does indeed give a solution of~\eqref{eq:pythtrip}.
The solutions under~\eqref{eq:rep2} can be written as follows:
\bea
&& a=2^j\,\alpha^2-\beta^2,\quad b = 2^{\frac{j+2}{2}}\,\alpha\beta,\quad c=2^j\,\alpha^2+\beta^2, \label{eq:rep2'} \\
&& a=\alpha^2-2^j\,\beta^2,\quad b = 2^{\frac{j+2}{2}}\,\alpha\beta,\quad c=\alpha^2+2^j\,\beta^2, \label{eq:rep2''}
\enea
with $j\ge 2$ even. To check uniqueness, consider first $b$~odd. Then only~\eqref{eq:rep1} applies, and $b,c$
together with the sign of~$a$ determine $\alpha^2,\beta^2$  from a quadratic equation. Since $\alpha,\beta>0$,
they are uniquely determined. If $b$~even, then \eqref{eq:rep2'} or~\eqref{eq:rep2''} apply. In fact, $j$ is clearly determined uniquely, and depending on whether $4|a+c$ or $4|c-a$ exactly one of the representations~\eqref{eq:rep2'} or~\eqref{eq:rep2''} holds (note that one cannot have both $4|a+c$ and $4|c-a$). It is now clear that~$\alpha,\beta$ are unique. 
\eproof

\noindent
We now use this lemma to estimate the number of solutions of~\eqref{eq:tang}.

\begin{lemma}
\label{lem:cardsol}
The number of integer solutions of~\eqref{eq:tang} is no larger than~$k_1\,N\log N$, and no smaller than
$k_2 N\log N$, where $k_1,k_2$ are absolute multiplicative constants. 
Under the additional restriction that $N/2<c$ the number
of solutions to~\eqref{eq:tang} is at least $k_2 N$. 
\end{lemma}
\noindent
\underline{Proof:} We first deal with the upper bound. Notice that by symmetry it suffices to
bound the number of solutions of~\eqref{eq:tang} for which $b,c>0$. 
Consider first the number of solutions of~\eqref{eq:tang} under the additional restriction
\begin{equation}
\label{eq:addrest}
\gcd(|a|,b,c)=1.
\end{equation}
In that case, Lemma~\ref{lem:pythago} applies. Since $1\le c\le N$, 
the representation~\eqref{eq:rep1} requires that $1\le \alpha,\beta\le \sqrt{2N}$. Hence, 
\eqref{eq:rep1} cannot generate more than~$2N$ many solutions. Similar considerations show that
\eqref{eq:rep2'} and~\eqref{eq:rep2''} cannot produce more than~$N$ solutions each. 
Indeed, by symmetry it suffices to treat~\eqref{eq:rep2'}. Then
\[ 1\le\beta \le \sqrt{N}, \quad 1\le \alpha\le 2^{-\frac{j}{2}} \sqrt{N}\]
so that the contribution from~\eqref{eq:rep2'} does not exceed
\[ \sum_{\substack{ j\ge 2\\ j \text{\ even}}} 2^{-\frac{j}{2}}\,N \le N.\]
To remove~\eqref{eq:addrest}, consider all solutions of~\eqref{eq:tang}, $b,c>0$,  with
\[ \gcd(|a|,b,c)=D \text{\ \ where\ \ } 1\le D\le N.\]
Then by the previous case $D=1$, one concludes that there are no more than
\[ \les \sum_{D=1}^N \frac{N}{D} \les N\log N\]
many solutions of~\eqref{eq:tang}, as desired. \\
For the lower bound, we need to count how many pairs $(\alpha,\beta)\in[1,\sqrt{N}]^2$ there are
with $\alpha,\beta$ odd and relatively prime. But it is a well-known property ot the Euler $\phi$ function
that there are $\gtrsim N$ such pairs. Indeed, for any positive integer $M$, 
\begin{equation}
\label{eq:euler}
\sum_{k=1}^M \phi(2(2k+1)) \le \#\Bigl\{(n,m)\in [1,2(2M+1)]^2 \:|\: \gcd(n,m)=1,\;2\not|n,\;2\not|m\Bigr\}.
\end{equation}
The same technique that is used to show that $\sum_{m=1}^M \phi(m)=\frac{3M^2}{\pi^2}+O(M\log M)$,
see Section~18.5 in~\cite{HW},  
yields a lower bound $\gtrsim M^2$ for the left-hand side, as desired. 
This shows that the number of solutions of~\eqref{eq:tang} 
under the additional restriction~\eqref{eq:addrest} is at least $\gtrsim N$. 
Hence, summing over $1\le D\le N$ where
$D=\gcd(|a|,b,c)$, leads to the lower bound
\[ \sum_{D=1}^{N/2} N/D \gtrsim N \log  N,\]
as claimed. The final claim of the lemma follows by taking $\alpha,\beta \in[\sqrt{N}/2,\sqrt{N}]$  odd and
relatively prime.  One checks that the same arguments based on~\eqref{eq:euler} as before apply and we are done.
\eproof

\noindent 
Define a $\mu$-fold point to be a point at which between $\mu$ and~$10\mu$ circles are tangent.
In what follows we will count $\mu$-fold points together with their multiplicity, i.e., if
$k$ lines $\ell_1,\ldots,\ell_k$ meet at a point $Q$ and each line contains between $\mu$ and~$10\mu$
points that are the centers of circles that have $Q$ as common point of tangency, then we count $Q$ 
as $k$ points.\\
Lemma~\ref{lem:cardsol} now yields the following.

\begin{lemma}
\label{lem:incid}
Let $\delta>0$ be small and consider the family of circles 
\begin{equation}
\label{eq:circfam}
 \C:=\{C(j\delta,k\delta,\ell\delta)\:|\: 1\le j,k\le N,\; N/2\le \ell \le N\} 
\end{equation}
where $N=[\delta^{-1}]$. Then 
\bea
\label{eq:tangpairs} 
&& |\C|^{\frac43}\log |\C| \sim \card\{ (C,C')\in\C^2\:|\: C,C'\text{\ \ are tangent}\}, \\
&& |\C|^{\frac43}\log |\C| \les \card\{ (C,C')\in\C^2\:|\: C,C'\text{\ \ are tangent},\quad d(C,C')\sim1\}.
\nn
\enea 
Moreover,  the number of $1$-fold points for the family~$\C$ is~$\gtrsim |\C|^{\frac43}$. 
\end{lemma}
\noindent 
\underline{Proof:} Rescaling by $1/\delta$ yields the family of circles with integer lattice points
in $[1,N]^2$ as centers and integer radii between $N/2$ and~$N$ considered above. 
We now work with this rescaled family.
Given $C=C(x,y,r)\in\C$ arbitrary, it follows that the number of circles~$C'\in\C$, $C'\not=C$ which are
 tangent to~$C$ is larger than the number of integer solutions $(a,b,c)$ of
\begin{equation}
\nn 
a^2+b^2=c^2\not=0,\quad |a|,|b|,|c|\le \frac{N}{2}
\end{equation}
and no larger than the number integer solutions $(a,b,c)$ of
\begin{equation}
\nn 
a^2+b^2=c^2\not=0,\quad |a|,|b|,|c|\le N,
\end{equation}
see \eqref{eq:tang}.
Since $|\C|\sim N^3$, the estimates~\eqref{eq:tangpairs} follows from Lemma~\ref{lem:cardsol}.
By definition, a $\mu$-fold point is the same as 
\begin{equation}
\label{eq:clamshell}
\{C_j=C(x_j,y_j,r_j)\in\C\:|\: 1\le j\le M,\quad C_1,\ldots,C_M \text{\ tangent at one point}\}
\end{equation}
where $\mu\le M\le 10\mu$, and $M$ is maximal (in the sense that one cannot enlarge the 
set in~\eqref{eq:clamshell}). 
It is therefore enough to show that the number of $M$-tuples as in~\eqref{eq:clamshell} with $M\sim1$ is
at least $|\C|^{\frac43}$. For this it suffices to show that the number of solutions to
\[ a^2+b^2=c^2\not=0,\quad N/4\le c \le N/2,\quad \gcd(|a|,b,c)=1 \]
is $\gtrsim N$ (note that if $C(x_j,y_j,r_j)$ are tangent to $C(x,y,r)$ at a common point for $1\le j\le J$, then
$(x_j-x,y_j-y,r_j-r)$ are linearly dependent for $1\le j\le J$). 
 But that is precisely the final assertion of Lemma~\ref{lem:cardsol}, 
so we are done.\eproof

\noindent
The example from Lemma~\ref{lem:incid} does not lend itself to ``fattening'' up the circles in any reasonable way. More precisely, it is easy to see that
\begin{equation}
\label{eq:deltapairs} 
|\C|^{\frac53}  \les \card\Bigl\{ (C,C')\in\C^2\:|\: \Delta(C,C')<\frac{\delta}{100},\quad d(C,C')\sim1\Bigr\},
\end{equation} 
and the same holds for any other small absolute constant instead of $\frac{1}{100}$. 

It is, however, a simple matter to produce a random example with the desired properties (and without using any of the arithmetic considerations from before): Let $\C_0$ be a collection of $\delta$-separated circles $\{C(x_j,y_j,r_j)\}_{j=1}^N$ which is maximal, i.e., the points
$(x_j,y_j,r_j)\in [0,1]^2\times [1/2,1]$ form a $\delta$-net. Clearly, $N\sim\delta^{-3}$. Moreover, each circle
in this family is $10\delta$-tangent to about~$\delta^{-2}$ others, and the multiplicity of each $\delta$-rectangle (of which there are about $\delta^{-2}$ many) is about~$\delta^{-\frac32}$. Hence the total number of $\delta$-incidences is~$\sim N^{\frac53}$. The idea is now to choose each
circle with probability $p=A\delta^{\frac32}$, with $A$ some large absolute constant. 
To obtain a bipartite situation, we break up $\C_0$ into two pieces at a distance~$\sim1$ from one another. 
This leads to two random sets $\W$ and~$\B$ of circles at a mutual distance~$\sim1$. Choosing $\delta$ sufficiently
small one obtains that
\begin{equation}
\label{eq:badevent}
\Prob\Bigl[ |\W| < \delta^{-\frac32} \text{\ \ or\ \ }|\B| < \delta^{-\frac32} \Bigr] < e^{-\delta^{-1}},
\end{equation}
see Lemma~\ref{lem:bern} above. 
Denote the complement of the event in~\eqref{eq:badevent} by~$\good_0$. Define a $\delta$-rectangle to 
be~{\em good} provided it has multiplicity $(\ge1,\ge1)$, but not $(\ge A^2,\ge1)$ or $(\ge1,\ge A^2)$. 
The number of circles tangent to a given $\delta$-rectangle is Poissonian with mean 
\[\sim p\delta^{-\frac32}\sim A.\]
Thus for sufficiently large $A$ the probability that a given rectangle is good is at least~$\half$. 
It follows that the conditional expectation of good rectangles relative to the event $\good_0$ satisfies
\begin{equation}
\label{eq:bedingt} 
\Erw[\#\;\text{good rectangles}\:|\:\good_0] \gtrsim \delta^{-2}.
\end{equation}
This already provides an example for the sharpness of the $\les |\C|^{\frac43}$ bound on the number of  $(\delta,1)$-rectangles of multiplicity $(\ge1,\ge1)$ in the bipartite setting.  Indeed, $\delta^{-2}=(\delta^{-\frac32})^{\frac43}$. One can also achieve in addition that the total number of incidences does not exceed~$|\C|^{\frac43+\epsilon}$. To see this, let $A$ be fixed and set
\[ T_{ij} = \chi_{[C_i,\,C_j\text{\ are $5\delta$-tangent}]}\]
for every $i\not= j$. Then $\Erw\; \sum_{i\not=j}\,T_{ij}\sim p^2\,N^{\frac53}$, and one also checks that
\[ \Erw\Bigl(\sum_{i\not=j} T_{ij}\Bigr)^s \le C_s\,N^{s\frac53}\,p^{2s}.\]
Hence
\[
\Prob\Bigl[ \sum_{i\not=j} T_{ij} > K\,N^{\frac53}\,p^2 \Bigr] \le C_s\, K^{-s}.
\]
Let the complement of the event on the left-hand side be $\good_1$
where we have set $K=\delta^{-\epsilon}$ for some arbitrary but fixed $\epsilon>0$. Then one has
\[ \Prob[\good_1^c] \le C_s\delta^{s\epsilon},\]
and therefore also
\[ \Erw[\#\;\text{good rectangles}\:|\:\good_0\cap\good_1] \gtrsim \delta^{-2}\]
provided $s\epsilon>3$, say,  and $\delta$ sufficiently small. 
Thus, (with positive probability) there exist $\W$ and $\B$ so that $|\W|\sim|\B|\sim\delta^{-\frac32}$, 
\[ \#\;\text{incidences} = \sum_{i\not=j} T_{ij} < \delta^{-2-\epsilon} \text{\ \ \ and\ \ \ } \#\;\text{good rectangles} \gtrsim \delta^{-2}. \]

\medskip\noindent
\underline{Acknowledgment:} The author is indebted to Alex Iosevich, Igor Rodnianski, and Andreas Seeger
for pointing out several references.

\medskip\noindent
\underline{Author's address:} Department of Mathematics, 253-37 Caltech,\\
Pasadena, CA 91125, U.S.A. \\
email: schlag@its.caltech.edu

\end{document}